\documentclass{amsart}
\usepackage{amssymb}
\usepackage{amsfonts}
\usepackage{fancyhdr}
\usepackage{hyperref}
\usepackage{amsmath}

\usepackage{bbm}
\usepackage{graphicx}
\usepackage{float}
\newtheorem{lemma}{Lemma}

\newtheorem{theorem}{Theorem}
\theoremstyle{plain}

\newtheorem{corollary}{Corollary}

\newtheorem{definition}{Definition}

\newtheorem{proposition}{Proposition}
\newtheorem{remark}{Remark}

\numberwithin{equation}{section}

\newcommand\ontopof[2]{\genfrac{}{}{0pt}{}{#1}{#2}}

\title{  Improved convergence rates for some kernel random forest algorithms }
\author{Iakovidis Isidoros}
\address{Dipartimento di Matematica, Universit\`a di Bologna, 40126, Bologna, Italy}
\email{isidoros.iakovidis2@unibo.it}

\author{Nicola Arcozzi}
\address{Dipartimento di Matematica, Universit\`a di Bologna, 40126, Bologna, Italy}
\email{nicola.arcozzi@unibo.it}

\thanks{The first author was funded by the PhD scholarship PON-DOT1303154-3, the second author was partyally funded by project GNAMPA 2022 “HOLOMORPHIC FUNCTIONS IN ONE AND SEVERAL VARIABLES” and by Horizon H2020-MSCA-RISE-2017 project 777822 "GHAIA"}

\date{October 2023}

\subjclass[2020]{Primary: 	62G05 ; Secondary: 62G20, 46E22, 43A75}
\keywords{Nonparametric statistics, Kernel random forests, reproducing kernel spaces}

\begin{document}

\maketitle
\begin{abstract}
 Random forests are notable learning algorithms first introduced by Breinman in 2001, they are widely used for classification and regression tasks and their mathematical properties are under ongoing research. We consider a specific class of random forest algorithms related to kernel methods, the so-called  KeRF (Kernel Random Forests.) In particular, we investigate thoroughly two explicit algorithms, designed independently of the data set, the centered KeRF and the uniform KeRF. In the present article, we provide an improvement in the rate of convergence for both algorithms and we explore the related reproducing kernel Hilbert space defined by the explicit kernel of the centered random forest.
\end{abstract}

\section{Introduction}
Random forests are a class of non-parametric statistic machine learning algorithms used for regression and classification tasks. Random forest algorithms have the capability to perform with high accuracy in high dimensional sparse tasks avoiding overfitting. In particular random forests are considered to be among the most accurate learning algorithm classes for general tasks. They are routinely used in many fields including bio-informatics \cite{Bioinf}, economics \cite{Econ}, biology \cite{Biol}, and linguistics \cite{Ling}. 

The most widely used random forest algorithm was introduced by Breinman \cite{BR}, who was inspired by the work on random subspaces of Ho \cite{Ho},  and geometrical feature selection of Amit and Geman \cite{Gem}. In Breinma's random forest, the trees are grown based on the CART procedure(Classification And Regression Trees)  where both splitting directions and training sets are randomized. However, despite the few parameters that need to be tuned \cite{Li}, \cite{Gen}, their mathematical properties are still areas of active research \cite{Lin}, \cite{Dev}. A significant distinction among the class of random forest algorithms consists in the way each individual tree is constructed, and, in particular, the dependence of each tree on the data set. Some of the researchers consider random forests designed independently from the data set\cite{Lug},\cite{As},\cite{Den}. \\
In 2012 Biau in \cite{Bia} studied a random forest model proposed by Breinman, where the construction is independent of the data set, called in literature {\it centered random forest}. In \cite{Bia} an upper bound on the rate of consistency of the algorithm and its adaption to sparsity were proven. More precisely, about the first item, for a data set of $n$ samples in a space of dimension $d$, the estimate for the convergence rate was $ \mathcal{O}\left( n^{-\frac{1}{d \frac{4}{3}\log{2}+1}}\right) $. In 2021 Klusowski in \cite{Sharp_K} improved the rate of convergence to $ \mathcal{O}\left( (n\log^{\frac{d-1}{2}} {n})^{ -( \frac{1+\delta}{d\log{2}+1})}  \right),  $  where $\delta$ is a positive constant that depends on the dimension of the feature space $d$ and converges to zero as $d$ approaches infinity. In addition, in the same paper, Klusowski proved that the rate of convergence of the algorithm is sharp, although it fails to reach the minimax rate of consistency over the class of the Lipschitz functions \cite{minimax} $ \mathcal{O}\left(n^{\frac{-2}{d+2}} \right) .$ 
There is also important work on the consistency of algorithms that depend on data \cite{Ment}, \cite{Wag},\cite{CONS}.
For a comprehensive overview of both theoretical and practical aspects of the random forests see e.g. \cite{BS}, which surveys the subject up to 2016.

An important tool for algorithmically manipulating random forests is by means of kernel methods. Already Breinman \cite{KerB} pointed this out, and later this was formalized by Geurts et al. in \cite{Geu}. In the same direction Scornet in \cite{S} defined KeRF (Kernel Random Forest) by modifying the original algorithm, and providing theoretical and practical results. In particular, in his important work, Scornet provided explicit kernels for some generalizations of algorithms, their rate of consistency, and comparisons with the corresponding random forests. Furthermore, in the very recent \cite{Interp} Arnould et al. investigated the trade-off between interpolation of several random forest algorithms and their consistency results. Moreover, it is proven that the centered KeRF is consistent when the trees are grown in full size when the feature set is uniformly distributed.

In the first part of the paper, we provide the notation and the definitions of the centered and uniform random forests and their corresponding kernel-related formulations. In addition, we improve the rate of consistency for the centered KeRF algorithm. Let $k\ge 1$ be the depth of the trees used to estimate
the target variable $Y$ (see Section \ref{SecPetronio} for definitions and notation).
\begin{theorem}\label{TheoPetronio} Suppose that $X=(X_1,\dots,X_d)$ and $Y$ are related by
  $   \textbf{Y}=m(\textbf{X}) +\epsilon$ where: $\epsilon$ is a zero mean Gaussian noise with finite variance independent of $\textbf{X} $, $\textbf{X} $ is uniformly distributed in $ [0,1]^d$, and $m$ is a regression function, which we assume to be Lipschitz. 
  As $k\to \infty$, there exists a constant $ \Tilde{C}$ such that, for every $n>1$ and for every $x \in [0,1]^d$, 
\[  \mathbb{E}(\tilde{m}^{Cen}_{\infty,n}(x)- m(x) )^2 \leq   \Tilde{C} n^{-\bigl(\frac{1}{1+d \log{2}}\bigr)} (\log{n}) .\]
\end{theorem}
Here, $m(x)=\mathbb{E}\left[Y|X=x\right]$ is the predicted value of $Y$ for $X=x\in [0,1]^d$, while $\tilde{m}^{Cen}_{\infty,n}(x)$ is the estimate for $m$ provided by the kernel version of the centered random forest algorithm, as $k\to\infty$.

Similarly, with $\tilde{m}^{Un}_{\infty,n}(x)$ playing for the uniform KeRF algorithm the role $\tilde{m}^{Cen}_{\infty,n}(x)$ had above, we have:
\begin{theorem}\label{TheoAmbrogio} Let $\textbf{X}$, $\textbf{Y}$, $m$, $\epsilon$, and $k$ be as in Theorem \ref{TheoPetronio}, with
  $   \textbf{Y}=m(\textbf{X}) +\epsilon$. As $k\to \infty$, there exists a constant $ \Tilde{C}$ such that for every $n>1$ , for every $x \in [0,1]^d$ 
\[  \mathbb{E}(\tilde{m}^{Un}_{\infty,n}(x)- m(x) )^2 \leq   \Tilde{C} n^{-\bigl(\frac{1}{1+\frac{3}{2}d \log{2}}\bigr)} (\log{n}) .\]
\end{theorem}
Moreover, in the next section, we provide numerical examples and experiments concerning the tuning parameter $k,$ which is the tree depth of the two kernel-based random forest algorithms, by comparing the $L_{2}$-error for different values and under specific assumptions on the data set.

In the final part of the article, we consider the reproducing kernel $K$ used in the centered KeRF algorithm {\it per se}. It is rewarding looking at it as defined on the finite Abelian group $\mathbb{Z}_{2}^{kd}$, where, as above, $d$ is the dimension of the vector $X$ and $k$ is the depth of the tree. 
By using elementary Fourier analysis on groups, we obtain several equivalent expressions for $K$ and its group transform, we characterize the functions belonging to the corresponding Reproducing Kernel Hilbert Space (RKHS) $H_K$, we derive results on 
multipliers, and we obtain bounds for the dimension of $H_K$, which is much smaller than what one might expect.

\section{Notation.}\label{SecPetronio}
A usual problem in machine learning is, based on $n$ observations of a random vector $(X,Y) \in \mathcal{X}\times\mathbb{R} \subseteq \mathbb{R}^{d}\times\mathbb{R}$, to estimate the function $m(x)=\mathbb{E}(Y | X=x).$ In classification problems, $Y$ ranges over a finite set.
In particular we assume that we are given a training sample $\mathcal{D}_{n}= \{(X_{1}, Y_{1}),...,(X_{n}, Y_{n})  \} $ of independent random variables, where $X_{i} \in [0,1]^{d} $ for every $i=1,...,n $ and $ Y \in \mathbb{R} $ with a shared joint distribution $\mathbb{P}_{X, Y}. $ The goal is using the data set to construct an estimate $m_{n}: \mathcal{X} \subseteq [0,1]^d \to \mathbb{R} $ of the function $m.$ Some {\it a priori} assumptions on the function $m$ ha to be made. Following \cite{S}, we suppose that $m$ belongs to the class of $L$ Lipschitz functions,
\[
|m(x)-m(x')|\le L \cdot \|x-x'\|.
\]
Here, as is \cite{S}, we consider on $\mathbb{R}^d$ the distance $ \|x-x'\|=\sum_{j=1}^d|x_j-x_j'|$.
\subsection{The Random Forest Algorithm.}
We provide here a description of the centered random forest. First, we define a {\it Random Tree} $\Theta$.
Start with a random variable $\Theta^0$, uniformly distributed in $\{1,\dots,d\}$, and split $I:=[0,1]^d=I^{\Theta^0}_0\cup I^{\Theta^0}_1$, where $I^{\Theta_0}_l=[0,1]\times\dots\times[l/2,(l+1)/2]\times\dots[0,1]$, where for $l=0,1$ the splitting was performed in the $\Theta^0$-th coordinate. Choose then random variables $\Theta^{1,l}$ ($l=0,1$), distributed as $\Theta^0$, and split
each $I^{\Theta^0}_l=I^{\Theta^0,\Theta^1}_{l,0}\cup I^{\Theta^0,\Theta^1}_{l,1}$, where, as before, the splitting is performed at the $\Theta_1$-th coordinate, and $I^{\Theta^0,\Theta^1}_{l,0}$ is the lower half of $I^{\Theta^0}_l$. Iterate the same procedure $k$ times.
In order to do that, we need random variables $\Theta^{j;\eta_0,\dots,\eta_j}$, with $\eta_l\in\{1,\dots,d\}$ and $j=1,\dots,k$. We assume that all such random variables are independent. It is useful think of ${\bf \Theta}=\{\Theta^{j;\eta_0,\dots,\eta_j}\}$  as indexed by a $d$-adic tree, and, in fact, we refer to ${\bf \Theta}$ as a random tree in $[0,1]^d$. We call {\it cells}, or {\it leaves}, each of the $2^k$ rectangles into which $[0,1]^d$ is split at the end of the $k^{th}$ subdivision.

Instead of splitting the $\Theta^l$-th factor into two equal parts, we might split it according to a random variable, or to the data set we are using. In these cases, we also refer to ${\bf \Theta}$ as a random tree, with a different notion of randomness and cells of different sizes.

A {\it random forest} is a finite collection of independent, finite random trees. More formally, a random forest is a collection $\{{\bf \Theta}_1,\dots,{\bf \Theta}_M\}$ of random trees.
Next, we define a random tree in a more general context than the centered random forest. 
Let's assume $\Theta_{i} $ for $i=1,..., M$ is a collection of independent random variables, independent of the data set $\mathcal{D}_{n} $ distributed as $\Theta.$ The random variables $\Theta_{i} $ correspond to sample the training set or select the positions for splitting.

  \begin{figure}[H]
    \centering
    \includegraphics[width=0.4\textwidth]{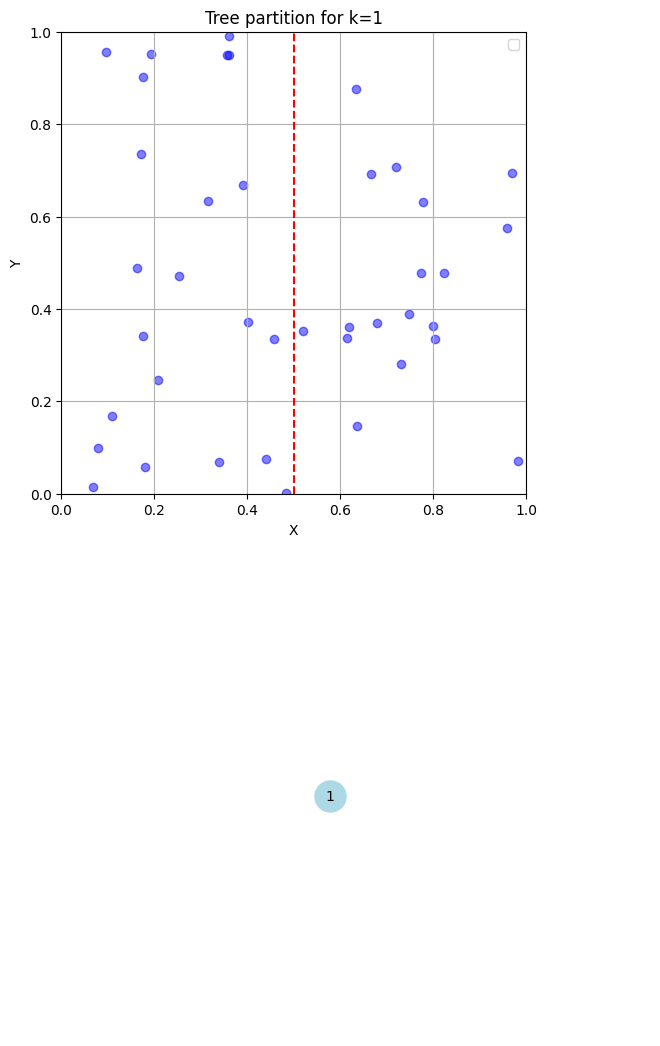}
             \caption{Centered algorithm with tree level  $k=1 $ with the convention that $1$ corresponds to $x$ axis and $2$ to the $y$ axis.}
    \label{fig:plot_label}
\end{figure}

  \begin{figure}[H]
    \centering
    \includegraphics[width=0.4\textwidth]{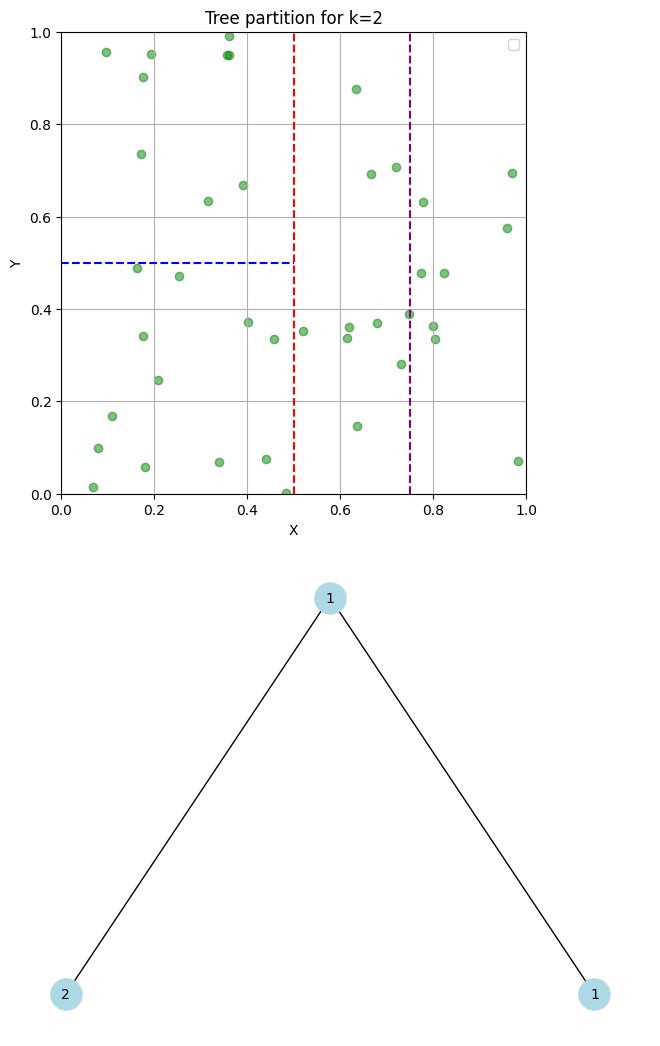}
     \caption{Centered algorithm with tree level  $k=1 $ with the convention that $1$ corresponds to $x$ axis and $2$ to the $y$ axis.}
    \label{fig:plot_label}
\end{figure}

\begin{definition}
For the $j$-th tree in the family of forests, the predicted value x will be denoted by 
\[ m_{n,\Theta_{j},\mathcal{D}_{n}}(x)=\sum_{i =1}^{n} \frac{\mathbbm{1}_{X_{i}\in A_{n,\Theta_{j},\mathcal{D}_{n} }( x)} Y_{i}}{N_{n,\Theta_{j},\mathcal{D}_{n}}( x) } .\]

Where $ A_{n,\Theta_{j},\mathcal{D}_{n}}( x) $ is the cell containing $x$ and $N_{n,\Theta_{j},\mathcal{D}_{n}}( x)$
is the number of points that fall into the cell that $x$ belongs to.
\end{definition}
Only one term in the sum is non-vanishing, and it is the average of the measured values of $Y$ in the cell containing $x$, which is, this is the hope, a good guess for the target value corresponding to $x$.

Now we can define the forest for the predicted variable $x$ as the average  of $M-$ random trees:
\begin{definition}
The finite $M$ forest is 
\[ m_{M,n}(x)= \frac{1}{M}\sum_{j=1}^{M} m_{n,\Theta_{j} ,\mathcal{D}_{n}}(x) . \]
\end{definition}
From a modeling point of view, we let $M\to \infty$ and  consider the infinite forest estimate 
\[m_{\infty,n,\mathcal{D}_{n}}(x  )=\mathbb{E}_{\Theta}( m_{n,\Theta,\mathcal{D}_{n} }(x)) . \]
The convergence holds almost surely by the law of the large numbers (Breinman) \cite{AL.SUR}, (Scornet) \cite[Theorem 3.1]{AL.SUR2}.

\subsection{KeRF Random Forest algorithm}

In 2016, E.Scornet \cite{S} introduced  kernel methods in the random forest world (KeRF), producing a kernel-based algorithm, together with estimates on how this compares with the old one, described above.\\
To understand the intuition behind KeRF random forests, we reformulate the random forest algorithm.\\
For all $x \in [0,1]^d ,$
\[  m_{M,n}(x)= \frac{1}{M}\sum_{j=1}^{M} \big(\sum_{i=1}^n \frac{\mathbbm{1}_{X_{i}\in A_{n,\Theta_{j},\mathcal{D}_{n} }( x)} Y_{i}}{N_{n,\Theta_{j},\mathcal{D}_{n}}( x) }\big). \]
Therefore we can define the weights of every observation $Y_{i}$ as
\[W_{i,j,n}(x)=\frac{\mathbbm{1}_{X_{i}\in A_{n,\Theta_{j},\mathcal{D}_{n} }( x)}}{N_{n,\Theta_{j},\mathcal{D}_{n}}( x) } .\]
Hence it is clear that the value of weights changes significantly concerning the number of points in each cell. A way to overcome this nuisance is by simultaneously considering all tree cells containing $x$, as the tree is randomly picked in the forest.\\
For all $x \in [0,1]^{d},$
\[  \Tilde{m}_{M,n,\Theta_{1},\Theta_{2},...,\Theta_{M}}(x )=\frac{1}{\sum_{j=1}^{M}N_{n,\Theta_{j}}(x)} \sum_{j=1}^{M}\sum_{i=1}^{n}Y_{i}\mathbbm{1}_{X_{i}\in A_{n,\Theta_{j}}(x)} . \]
This way, empty cells do not affect the computation of the prediction function of the algorithm.

It is proven \cite{S}, that this representation has indeed a kernel representation.
\begin{proposition}[Scornet \cite{S}, Proposition 1] 
    For all $x \in [0,1]^{d}$ almost surely, it holds
    \[  \Tilde{m}_{M,n,\Theta_{1},\Theta_{2},...,\Theta_{M}}(x )= \frac{\sum_{i=1}^{n}K_{M,n}(x,X_{i})Y_{i}}{\sum_{i=1}^{n}K_{M,n}(x,X_{i})} . \]
    where 
    \[ K_{M,n}(x,z)= \frac{1}{M}\sum_{i=1}^{M}\mathbbm{1}_{ x \in A_{n,\Theta_{i},\mathcal{D}_{n}}( z  ) }. \]
    is the proximity function of the $M$ forest
\end{proposition}
The {\it infinite random forest} arises as the number of trees tends to infinity.
\begin{definition}
   The {\bf infinite KeRF} is defined as:
   \[ \Tilde{m}_{\infty,n}(x)= \lim_{M\to \infty}\Tilde{m}_{M,n}(x,\Theta_{1},\Theta_{2},...,\Theta_{M} ) .  \]
\end{definition}
The extension of the kernel follows also in the infinite random forest.

\begin{proposition}[Scornet \cite{S}, Proposition 2]
 Almost surely for all $x,y \in [0.1]^{d}$
\[\lim_{M\to \infty} K_{M,n}(x,y)=K_{n}(x,y), \]
 where
 \[ K_{n}(x,y)= \mathbb{P}_{\Theta} (x\in A_{n}(y,\Theta)),\]
 where the left-hand side is the probability that $x$ and $y$ belong to the same cell in the infinite forest.
\end{proposition}
\subsection{The Centred Random Forest/KeRF Centred Random Forest and the Uniform Random Forest / Kernel Random Forest.}

An estimate function $m_{n} $ of $m$ is consistent if $ L_{2}$ convergence holds,
\[\mathbb{E}( m_{n}(x) - m (x) )^{2} \to 0 ,\] as $n \to \infty.$

In the centered and uniform forest algorithms, the way the data set $\mathcal{D}_{n}$ is partitioned is independent of the data set itself.  
\subsubsection{The centered random Forest/KeRF centered random forest.}

The centered forest is designed as follows.
\begin{itemize}
 \item[1)]  Fix $k \in \mathbb{N}.$  
 \item[2)] At each node of each individual tree choose a coordinate uniformly from $\{ 1,2,..d\}.$
 \item[3)] Split the node at the midpoint of the interval of the selected coordinate.
\end{itemize}
Repeat step 2)-3) $k$ times. At the end, we have $2^{k}$ leaves, or cells.
Our estimation at a point $x$ is achieved by averaging the $Y_{i}$ corresponding to the $X_{i}$ in the cell containing $x.$\\
E.Scornet in \cite{S} introduced the corresponding kernel-based centered random forest providing explicitly the proximity kernel function.

\begin{proposition}
    A centered random forest kernel with $k \in \mathbb{N}$ parameter has the following multinomial expression \cite[Proposition 5]{S}.
    \[ K^{Cen}_{k}(x,z) = \sum_{\sum_{j=1}^{d}k_{j}=k } \frac{k!}{k_{1}!...k_{d}!}(\frac{1}{d})^k  \prod_{j=1}^{d}  \mathbbm{1}_{ \left \lceil{2^{k_{j}}x_{j}}\right \rceil =\left \lceil{2^{k_{j}}z_{j}}\right \rceil } . \]
\end{proposition}
Where $K^{Cen}_{k}$ is the Kernel of the corresponding centered random forest. 

\subsubsection{The uniform random forest / Kernel random forest.}
Uniform Random forest was introduced by Biau et al. \cite{Lug} and is another toy model of Breinman's random forest as a centered random forest. The algorithm forms a partition in $ [0,1]^{d}$ as follows:
\begin{itemize}
 \item[1)]  Fix $k \in \mathbb{N}.$  
 \item[2)] At each node of each individual tree choose a coordinate uniformly from $\{ 1,2,..d\}.$
 \item[3)] The splitting is performed uniformly on the side of the cell of the selected coordinate.
\end{itemize}
Repeat step 2)-3) $k$ times. At the end, we have $2^{k}$ leaves.
Our final estimation at a point $x$ is achieved by averaging the $Y_{i}$ corresponding to the $X_{i}$ in the cell $x.$\\
Again Scornet in \cite[Proposition 6]{S} proved the corresponding kernel-based uniform random forest. 
\begin{proposition}
    The corresponding proximity kernel for the uniform Kernel random forest for parameter $k \in \mathbb{N}$ and $x\in [0,1]^d $ has the following form:
    \[ K^{Un}_{k}(0,x)= \sum_{\sum_{j=1}^{d}k_{j}=k }\frac{k!}{k_{1}!...k_{d}!}(\frac{1}{d})^k \prod_{m=1}^{d} \bigg( 1-x_{m} \sum_{j=0}^{k_{m}-1}\frac{(-\ln{x_{m}})^{j}}{j!}  \bigg). \] with the convention that $\sum_{j=0}^{-1}\frac{(-\ln{x_{m}})^{j}}{j!}=0 $ and by continuity we can extend the kernel also for zero components of the vector.
\end{proposition}
Unfortunately, it is very hard to obtain a general formula for $K^{Un}(x,y)$ but we consider instead a translation invariant  KeRF uniform forest: 
\[ \title{m}^{Un}_{\infty,n}(x)= \frac{\sum_{i=1}^{n}Y_{i}K^{Un}_{k}(0,\lvert X_{i}-x \rvert ) }{\sum_{i=1}^{n}K^{Un}_{k}(0,\lvert X_{i}-x \rvert ) }.\]

   \section{Proofs of the main theorems.}
In this section, after providing some measure concentration type results \cite{meas_cons}, we improve the rate of consistency of the centered KeRF algorithm. The following lemmata will provide inequalities to derive upper bounds for averages of iid random variables.

\begin{lemma}\label{Bound_lemma}

  Let $X_{1},..., X_{n}$ be a sequence of independent and identically distributed random variables with $\mathbb{E}(X_{i})=0.$ Assuming also that there is a uniform bound for the $L_{1}$-norm and the supremum norm i.e. $\mathbb{E}(\rvert X_{i} \lvert) \leq C ,$ $\lvert \lvert X_{i}  \rvert \rvert_{\infty}  \leq CM$ for every $i=1,...,n.$ Then for every $t\in (0,1)$ 
  \[ \mathbb{P}\big( \{ \frac{\lvert\sum_{i=1}^{n}X_{i}  \rvert   }{n} \geq t\} \big)\leq 2 e^{- \Tilde{C}_{C}\frac{t^{2} n}{   M }} .\] for some positive constant $ \Tilde{C}_{C,M}$ that depends only on $C $ and $ M.$

  \begin{proof}
   $ \forall x \in [0,1]$ one has that $e^x \leq 1+x+x^2 .$ By using the hypothesis for every $\lambda \leq \frac{1}{CM} ,$

   \begin{align*}
       e^{\lambda X_{i}} &\leq 1 + \lambda X_{i} + (\lambda {X_i})^2 \quad \Rightarrow\\
       \mathbb{E}e^{\lambda X_{i}} &\leq 1  +\lambda^2 \mathbb{E}({X_i})^2\\
       &\leq 1 + \lambda^2 \lvert \lvert  X_{i}\rvert \rvert_{1} \lvert \lvert X_{i}\rvert \rvert_{\infty}\\
       &\leq 1 + \lambda^2 C^2 M\\
       &\leq e^{\lambda^2 C^2 M}.\\
   \end{align*}
   By the independence of the random variables $X_{i},$

   \begin{align*}
        \mathbb{E}e^{\sum_{i=1}^n \lambda X_{i}} &= \prod_{i=1}^{n}\mathbb{E}e^{\lambda X_{i}}\\
        &\leq \prod_{i=1}^{n}e^{\lambda^2 C^2 M}\\
        & = e^{n\lambda^2 C^2 M}.
   \end{align*}
   Therefore, by Markov inequality 
   \begin{align*}
       \mathbb{P}\big( \{ \frac{\sum_{i=1}^{n}X_{i}     }{n} \geq t\} \big) &\leq e^{-\lambda t n} \mathbb{E}e^{\sum_{i=1}^n \lambda X_{i}}\\
       & \leq e^{-\lambda t n} e^{n\lambda^2 C^2 M}\\
       &=e^{n\lambda^2 C^2 M-\lambda t n}.\\
   \end{align*}
   Finally if $C \geq \frac{1}{4} $ we choose, $\lambda=\frac{t}{2 C^2 M},$ otherwise for $ \lambda=\frac{t}{16 C M}$ 
   \[ \mathbb{P}\big( \{ \frac{\sum_{i=1}^{n}X_{i}     }{n} \geq t\} \big) \leq e^{- \tilde{C}_{C}\frac{t^{2} n  }{  M}}.  \]
   By replacing $X_{i}$ with $-X_{i}$ we conclude the proof.
  \end{proof}
\end{lemma}
 \begin{lemma}\label{Bounded_lemma2}
       Let $X_{1},..., X_{n}$ be a non-negative sequence of independent and identically distributed random variables with  $\mathbb{E}( X_{i} ) \leq 2,$ $\lvert \lvert X_{i}  \rvert \rvert_{\infty}  \leq M$ for every $i=1,...,n.$ Let also a sequence of independent random variables $\epsilon_{i}$ following normal distribution with zero mean and finite variance $\sigma^2,$ for every $i=1,...,n. $ We assume also that $\epsilon_{i}$ are independent from $X_{i} $  for every $i=1,...,n. $\\
       Then for every $t\in (0,1),$ 
 \[  \mathbb{P}\biggl(\frac{1}{n} \sum_{i=1}^{n} \lvert \epsilon_{i}X_{i}\rvert \geq t \biggr) \leq 2\exp{(-C t^{2}\frac{n}{M})}. \] with the positive constant $C$ depending only on $\sigma .$
\begin{proof}

\begin{align*}
    \mathbb{P}\biggl(\frac{1}{n} \sum_{i=1}^{n} \epsilon_{i}X_{i} \geq t \biggr)  &= \mathbb{P}\biggl( \exp\biggl(\frac{\lambda}{n} \sum_{i=1}^{n} \epsilon_{i}X_{i} \geq \exp(\lambda t) \biggr) \quad \text{for a positive } \lambda\\
    &\leq \exp(-\lambda t) \mathbb{E}\exp\biggl(\frac{\lambda}{n} \sum_{i=1}^{n}\epsilon_{i} X_{i}  \biggl)\quad \text{By Chebyshev's inequality} \\
    &= \exp(-\lambda t) \prod_{i=1}^{n}\mathbb{E}\exp\biggl(\frac{\lambda}{n} \epsilon_{i}X_{i} \biggl) \quad \text{By independence}\\
    &= \exp(-\lambda t)\prod_{i=1}^{n}\biggl( 1+ \sum_{k=2}^{\infty} \frac{\lambda^{k}\mathbb{E}X_{i}^{k}\mathbb{E}\epsilon_{i}^{k}}{n^{k}k!}\biggl) \\
    &\leq \exp(-\lambda t)\prod_{i=1}^{n}\biggl( 1+ \frac{2}{M}\sum_{k=2}^{\infty} \frac{\lambda^{k}M^{k}\mathbb{E}\epsilon_{i}^{k}}{n^{k}k!}\biggr)\\
    &= \exp(-\lambda t)\prod_{i=1}^{n}\biggl( 1+ \frac{2}{M} \biggr(  \mathbb{E}\exp\biggl(\frac{\lambda M}{n}\epsilon_{i}\biggr) -1\biggr) \biggr)\\
    &\leq \exp(-\lambda t)\prod_{i=1}^{n}\biggl( 1+ \frac{2}{M} \biggr( \exp\biggl(\frac{\lambda^2 \sigma^2 M^2}{n^2}\biggr) -1\biggr) \biggr)\\
    &= \exp(-\lambda t)\exp\biggl(\sum_{i=1}^{n}\biggl(\log\biggl( 1 + \frac{2}{M} \biggr( \exp\biggl(\frac{\lambda^2 \sigma^2 M^2}{n^2}\biggr) -1 \biggr)\biggr)\biggr)\biggr)\\
    &\leq \exp(-\lambda t) \exp\biggl( \sum_{i=1}^{n} \frac{2}{M}\biggl(\exp\biggl(\frac{\lambda^2 \sigma^2 M^2}{n^2}\biggr) -1\biggr)\biggr)\\
    &\leq \exp(-\lambda t) \exp\biggl( \frac{2n}{M}\biggl(\exp\biggl(\frac{\lambda^2 \sigma^2 M^2}{n^2}\biggr) -1\biggr)\biggr)\\
     &\leq \exp(-\lambda t) \exp\biggl( \frac{2n}{M}\biggl(2\frac{\lambda^2 \sigma^2 M^2}{n^2} \biggr)\biggr)\quad \text{for } \lambda\leq\frac{n}{ \sigma M}\\
     &=  \exp\biggl(-\lambda t+ \frac{4M}{n}\lambda^2 \sigma^2 \biggr).
\end{align*}
Finally we select $\lambda= \frac{tn}{8 M \sigma^2},$ when $ \sigma \geq \frac{1}{8}$ and $\lambda= \frac{tn}{M \sigma} ,$ when $ \sigma \leq \frac{1}{8} $
\[\mathbb{P}\biggl(\frac{1}{n} \sum_{i=1}^{n}\epsilon_{i}X_{i} \geq t \biggr)   \leq  \exp\biggl(- C \frac{t^2 n}{ M} \biggr).\]

Replacing $X_{i}$ with $-X_{i}$ we conclude the proof. 
\end{proof}
 \end{lemma}
\begin{theorem}
  $   \textbf{Y}=m(\textbf{X}) +\epsilon$ where $\epsilon$ is a zero mean Gaussian noise with finite variance independent of $\textbf{X} $. Assuming also that $\textbf{X} $ is uniformly distributed in $ [0,1]^d$ and $m$ is a Lipschitz function. Providing $k\to \infty$, there exists a constant $ \Tilde{C}$ such that for every $n>1$ , for every $x \in [0,1]^d$ 
\[  \mathbb{E}(\tilde{m}^{Cen}_{\infty,n}(x)- m(x) )^2 \leq   \Tilde{C} n^{-\bigl(\frac{1}{1+d \log{2}}\bigr)} (\log{n}) .\]

 \begin{proof}

Following the notation in \cite{S}, let $x \in [0,1]^d$, $ \lVert m \rVert_{\infty} = \sup_{x \in [0,1]^d} | m(x) |$,
and by the construction of the algorithm
\[ \tilde{m}_{n,\infty}^{Cen}(x) = \frac{\sum_{i=1}^{n}Y_{i}K_{k}(x,X_{i})}{\sum_{i=1}^{n}K_{k}(x,X_{i})} .\]
Let \[A_{n}(x) = \frac{1}{n}\sum_{i=1}^{n} \left( \frac{ Y_{i}K_{k}(x,X_{i}) - \mathbb{E}(Y K_{k}(x,X))} {\mathbb{E}(K_{k}(x,X))} \right),\]
\[B_{n}(x) = \frac{1}{n}\sum_{i=1}^{n} \left( \frac{K_{k}(x,X_{i}) - \mathbb{E}( K_{k}(x,X))} {\mathbb{E}(K_{k}(x,X))} \right),\]
and \[M_{n}(x) = \frac{\mathbb{E}(YK_{k}(x,X))}{\mathbb{E}( K_{k}(x,X))}.\]
Hence, we can reformulate the estimator as
\[ \tilde{m}_{n,\infty}^{Cen}(x) = \frac{M_{n}(x)+A_{n}(x)}{B_{n}(x)+1}. \]
Let $t \in (0,\frac{1}{2})$ and the event $C_{t}(x)$ where $\{A_{n}(x),B_{n}(x) \leq t\}.$
\begin{align*}
\mathbb{E}(\tilde{m}_{n,\infty}^{cc}(x) - m(x))^2 &= \mathbb{E}(\tilde{m}_{n,\infty}^{cc}(x) - m(x))^2 \mathbbm{1}_{C_{t}(x)} + \mathbb{E}(\tilde{m}_{n,\infty}^{cc}(x) - m(x))^2 \mathbbm{1}_{C_{t}^c(x)} \\
&\leq \mathbb{E}(\tilde{m}_{n,\infty}^{cc}(x) - m(x))^2 \mathbbm{1}_{C_{t}^c(x)} + c_{1} \left(1-\frac{1}{2d}\right)^{2k} + c_{2}t^2. \quad 
\end{align*}
Where the last inequality was obtained in \cite[p.1496]{S}
Moreover, in \cite{S},
\[ \mathbb{E}(\tilde{m}_{n,\infty}^{cc}(x) - m(x))^2 \mathbbm{1}_{C_{t}^c(x)} \leq c_{3}(\log{n})(\mathbb{P}(C_{t}^c(x)))^{\frac{1}{2}}.\]

In order to find the rate of consistency 
we need a bound for the probability      $ \mathbb{P}(C_{t}^c(x)). $
Obviously , \[  \mathbb{P}(C_{t}^c(x)) \leq    \mathbb{P}( \rvert A_{n}(x) \rvert >t)+  \mathbb{P}( \lvert B_{n}(x) \rvert >t) .\]
We will work separately to obtain an upper bound for both probabilities. 
\begin{proposition}\label{Prop_Bn}

Let $  \tilde{X}_{i} =   \frac{K_{k}(x,X_{i})}{\mathbb{E}(K_{k}(x,X))} -1$ a sequence of iid random variables. Then for any $t \in (0,1),$
\[ \mathbb{P}\big( \{ \frac{ \lvert \sum_{i=1}^{n}\tilde{X}_{i}     \rvert}{n} \geq t\} \big)= \mathbb{P}\big( \lvert B_{n}(x) \rvert \geq t \big) \leq 2e^{-\tilde{C}_{1} \frac{t^{2} n  }{2^k }}  \] for some positive constant $\tilde{C}_{1}.$
\begin{proof}
   It is easy to verify that $\mathbb{E} \tilde{X}_{i}=0, $ and 
   \[\lvert  \tilde{X}_{i} \rvert= \vert  \frac{K_{k}(x,X_{i})}{\mathbb{E}(K_{k}(x,X))} -1 \rvert \leq \frac{K_{k}(x,X_{i})}{\mathbb{E}(K_{k}(x,X))}+1,\] hence, $\mathbb{E} \lvert \tilde{X}_{i} \rvert \leq 2 . $\\
   Finally,
   \[\lvert \lvert  \tilde{X}_{i}  \rvert \rvert_{\infty}=\sup \{  \lvert \tilde{X}_{i} \rvert\}= \sup \{  \lvert \frac{K_{k}(x,X_{i})}{\mathbb{E}(K_{k}(x,X))} -1 \rvert\} \leq \frac{1}{\mathbb{E}(K_{k}(x,X))} \sup{ K_{k}(x,X_{i})}+1 \leq 2^k +1 \leq  2^{k+1}.\]
   By Lemma \ref{Bound_lemma},
   
   \[  \mathbb{P}\big( \{ \frac{ \lvert \sum_{i=1}^{n}\tilde{X}_{i}     \rvert}{n} \geq t\} \big)= \mathbb{P}\big( \lvert B_{n}(x) \rvert \geq t \big) \leq 2 e^{- \tilde{C}_{1} \frac{t^{2} n  }{ 2^k }}  . \]
\end{proof}
\end{proposition}
We need a bound for 
 $\mathbb{P}\big(\lvert A_{n}(x) \rvert > t \big)  $
 where,
 \[A_{n}(x)= \frac{1}{n}\sum_{i=1}^{n} \big( \frac{ Y_{i}K_{k}(x,X_{i})  - \mathbb{E}(Y K_{k}(x,X))} {\mathbb{E}(K_{k}(x,X))}  \big)  . \]

\begin{proposition}\label{Prop_An}
 Let  $\Tilde{Z_{i}}=\frac{ Y_{i}K_{k}(x,X_{i})  - \mathbb{E}(Y K_{k}(x,X))} {\mathbb{E}(K_{k}(x,X))}$  for $i=1,...,n $
    then for every $t \in (0,1),$
 \[  \mathbb{P}\big( \{ \frac{ \lvert \sum_{i=1}^{n}\tilde{Z}_{i}     \rvert}{n} \geq t\} \big)= \mathbb{P}\big( \lvert A_{n}(x) \rvert \geq t \big) \leq 4 e^{- C \frac{t^{2} n  }{2^k }},  \] for some constant $C$
 depending only on $\sigma, \lVert m \rVert_\infty .$ 
   \begin{proof}
   
\begin{align*}
    A_{n}(x) &= \frac{1}{n}\sum_{i=1}^{n} \bigg( \frac{ Y_{i}K_{k}(x,X_{i})  - \mathbb{E}(Y K_{k}(x,X))} {\mathbb{E}(K_{k}(x,X))}  \bigg)\\
    &=\frac{1}{n}\sum_{i=1}^{n} \bigg( \frac{ m(X_{i})K_{k}(x,X_{i})  - \mathbb{E}(m(X) K_{k}(x,X))} {\mathbb{E}(K_{k}(x,X))}  \bigg) +\frac{1}{n}\sum_{i=1}^{n} \bigg( \frac{ \epsilon_{i}K_{k}(x,X_{i})  - \mathbb{E}(\epsilon K_{k}(x,X))} {\mathbb{E}(K_{k}(x,X))}  \bigg)\\
    &=\frac{1}{n}\sum_{i=1}^{n} \bigg( \frac{ m(X_{i})K_{k}(x,X_{i})  - \mathbb{E}(m(X) K_{k}(x,X))} {\mathbb{E}(K_{k}(x,X))}  \bigg) +\frac{1}{n}\sum_{i=1}^{n} \bigg( \frac{ \epsilon_{i}K_{k}(x,X_{i}) }{\mathbb{E}(K_{k}(x,X))} \bigg) .
\end{align*}
    Therefore, 
    \begin{align*}
    \mathbb{P} \big( \lvert A_{n}(x) \rvert \geq t \big)  &\leq 
    \mathbb{P} \bigg( \bigg\lvert \frac{2}{n} \sum_{i=1}^{n} \frac{m(X_{i})K_{k}(x,X_{i}) - \mathbb{E}(m(X) K_{k}(x,X))}{\mathbb{E}(K_{k}(x,X))} \bigg\rvert \geq t \bigg) \\
    &\quad+ \mathbb{P} \bigg( \bigg\lvert \frac{2}{n} \sum_{i=1}^{n} \frac{\epsilon_{i}K_{k}(x,X_{i})}{\mathbb{E}(K_{k}(x,X))} \bigg\rvert \geq t \bigg).
\end{align*}
Let $ {Z_{i}}= \frac{ 2(m(X_{i})K_{k}(x,X_{i}) - \mathbb{E}(m(X) K_{k}(x,X)) )} {\mathbb{E}(K_{k}(x,X))} $ a sequence of iid random variables.
It is easy to verify that $\Tilde{Z_{i}} $ are centered and 
 \[\lvert  \tilde{Z}_{i} \rvert=    \lvert  \frac{ m(X_{i})K_{k}(x,X_{i}) - \mathbb{E}(m(X) K_{k}(x,X))} {\mathbb{E}(K_{k}(x,X))}  \rvert \leq 2 \lvert \lvert m \rvert \rvert_{\infty} \frac{K_{k}(x,X_{i})+\mathbb{E}(K_{k}(x,X))}{\mathbb{E}(K_{k}(x,X))} . \]\\
 Hence, 
 \[ \mathbb{E}\lvert  Z_{i} \rvert \leq 4   \lvert \lvert m \rvert \rvert_{\infty}\]
   Finally,
   \begin{align*}
       \lvert \lvert  Z_{i}  \rvert \rvert_{\infty} &=\sup \{  \lvert Z_{i} \rvert\}\\
       &= 2\sup \{ \lvert  \frac{ m(X_{i})K_{k}(x,X_{i}) - \mathbb{E}(m(X) K_{k}(x,X))} {\mathbb{E}(K_{k}(x,X))}  \rvert \}\\
    &\leq 2 \lvert \lvert m \rvert \rvert_{\infty}(2^k +1)\\
    &\leq  4 \lvert \lvert m \rvert \rvert_{\infty} 2^k .
   \end{align*}
 
   By Lemma \ref{Bound_lemma},
\[ \mathbb{P} \Bigg( \Bigg\lvert \frac{1}{n} \sum_{i=1}^{n} \frac{m(X_i)K_k(x,X_i) - \mathbb{E}(m(X)K_k(x,X))}{\mathbb{E}(K_k(x,X))} \Bigg\rvert \geq t \Bigg) \leq 2 e^{- C\frac{nt^2}{ 2^k}}.\]

Furthermore let $\Tilde{W}_{i} = \frac{2\epsilon_{i}K_{k}(x,X_{i})}{\mathbb{E}(K_{k}(x,X))}$ for $i=1,...,n $ a sequence of independent and identically distributed random variables. We can verify that for every for $i=1,...,n $:
   \[
\mathbb{E} \left(  \frac{2K_{k}(x,X_{i})}{\mathbb{E}(K_{k}(x,X))}  \right) \leq 2 .
\]
   Finally,
   \[
\sup \left\{ \left\lvert \frac{2 K_{k}(x,X_{i})}{\mathbb{E}(K_{k}(x,X))} \right\rvert \right\} \leq \frac{2}{\mathbb{E}(K_{k}(x,X))} \sup \{ K_{k}(x,X_{i}) \} \leq 2^{k+1} .
\]
By Lemma \ref{Bounded_lemma2} it is clear,

\[    \mathbb{P} \bigg( \bigg\lvert \frac{2}{n} \sum_{i=1}^{n} \frac{\epsilon_{i}K_{k}(x,X_{i})}{\mathbb{E}(K_{k}(x,X))} \bigg\rvert \geq t \bigg) \leq   2 e^{- C_{2}\frac{nt^2}{ 2^k}} .\]
We conclude the proposition by observing
\[   \mathbb{P} \big( \lvert A_{n}(x) \rvert \geq t \big) \leq 4 e^{- \min{ \{ C_{2},C\}}\frac{nt^2}{ 2^k}}   .\]

  \end{proof}

  \end{proposition}

Finally, let us compute the rate of consistency of the algorithm-centered KeRF.
By Propositions \ref{Prop_Bn},\ref{Prop_An} one has that
 \[  \bigg( \mathbb{P}(C_{t}^c(x))  \bigg)^{\frac{1}{2}}     \leq    \bigg( \mathbb{P}( \rvert A_{n}(x) \rvert >t)+  \mathbb{P}( \lvert B_{n}(x) \rvert >t)  \bigg)^{\frac{1}{2}} \leq c_{3} e^{- c_{4}\frac{nt^2}{ 2^k}},\] 
 for some constants $c_{3}, c_{4}$ independent of $k$ and $n.$

    Thus, 
   \[ \mathbb{E}( \tilde{m}_{\infty,n} -m(x) )^{2} \leq c_{1} \left(1-\frac{1}{2d}\right)^{2k} + c_{2}t^{2} + c_{3}\log{n} e^{-c_{4}t^{2}\frac{n}{2^k}}.\]

        We compute the minimum of the right-hand side of the inequality for $t \in (0,1),$ 
        \begin{align*}
            2c_{2}t - 2tc_{4}\log{n}c_{3} \frac{n}{2^k} e^{-c_{4}t^{2}\frac{n}{2^k}} &= 0 \quad \Rightarrow\\
            e^{-c_{4}t^{2}\frac{n}{2^k}} &= \frac{c_{2}}{c_{3}c_{4}} \frac{2^k}{n\log{n}} \quad \text{and}\\
            t^{2} &= \frac{1}{c_{4}} \frac{2^k}{n}\log{\left(\frac{c_{2}}{c_{3}c_{4}} \frac{n \log{n}}{2^k}\right)}.
        \end{align*}
        
        Hence, the inequality becomes

        \begin{align*}
             \mathbb{E}( \tilde{m}_{\infty,n} -m(x) )^{2} &\leq c_{1} \left(1-\frac{1}{2d}\right)^{2k} + c_{2} \frac{1}{c_{4}} \frac{2^k}{n}\log{\left(\frac{c_{2}}{c_{3}c_{4}} \frac{n \log{n}}{2^k}\right)} + c_{3} \log{n}\frac{c_{2}}{c_{3}c_{4}}\frac{2^k}{n\log{n}} \\
             &=c_{1} \left(1-\frac{1}{2d}\right)^{2k} + c_{2} \frac{1}{c_{4}} \frac{2^k}{n}\log{\left(\frac{c_{2}}{c_{3}c_{4}} \frac{n \log{n}}{2^k} e^{\frac{c_{2}}{c_{4}}}\right)}.
        \end{align*}
        
    For every $\epsilon_n \in (0,2] $ it holds, $\log{x} \leq \frac{1}{\epsilon_n} x^{\epsilon_n}.$ Then one has that

    \[\mathbb{E}( \tilde{m}_{\infty,n} -m(x) )^{2} \leq c_{1} \left(1-\frac{1}{2d}\right)^{2k} +  \frac{c_2 (   e^{\frac{c_{2}}{c_{4}}}  \frac{c_{2}}{c_{3}c_{4}} )^{n} }{c_4 \epsilon_n}\left(\frac{2^k}{n} (\log{n})^{\frac{\epsilon_n}{1-\epsilon_n}}\right)^{1-\epsilon_n}.\]

    We pick $k =  c(d) \log_{2}{\frac{n}{(\log{n})^{\frac{\epsilon_n}{1-\epsilon_n}}}},$
    
    thus, 

\[     \frac{c_2 (   e^{\frac{c_{2}}{c_{4}}}  \frac{c_{2}}{c_{3}c_{4}} )^{n} }{c_4 \epsilon_n}\left(\frac{2^k}{n} (\log{n})^{\frac{\epsilon_n}{1-\epsilon_n}}\right)^{1-\epsilon_n} \leq \frac{c'}{\epsilon_n} n^{(c(d)-1)(1-\epsilon_n)} \log{n}^{\epsilon_n(1-c(d))}, \]

for a constant $ c'$ independent of $n$ and,

\begin{align*}
 c_{1} \left(1-\frac{1}{2d}\right)^{2k}&=c_{1} \left(1-\frac{1}{2d}\right)^{2 ( c(d) \log_{2}{\frac{n}{(\log{n})^{\frac{\epsilon_n}{1-\epsilon_n}}}}) } \\
 &= c_{1} 2^{2c(d) \log_{2}{\left(1-\frac{1}{2d}\right)}\log_{2}{\frac{n}{(\log{n})^{\frac{\epsilon_n}{1-\epsilon_n}}}} }\\
 &= c_{1} n^{2c(d)  \log_{2}{\left(1-\frac{1}{2d}\right)}} \frac{1}{(\log{n})^{c(d) \frac{2\epsilon_n}{1-\epsilon_n}\log_{2}{\left(1-\frac{1}{2d}\right)}} }.\\
\end{align*}

Therefore, 
\[ c(d)= \frac{\epsilon_n-1}{2\log_{2}{\left(1-\frac{1}{2d}\right)}-(1-\epsilon_n)}.\]
 Finally,
\begin{align*}
     c_{1}n^{2c(d) \log_{2}{\left(1-\frac{1}{2d}\right)}} \frac{1}{(\log{n})^{c(d) \frac{2\epsilon_n}{1-\epsilon_n}\log_{2}{\left(1-\frac{1}{2d}\right)}}}
&= c_{1} n^{\frac{2(\epsilon_n-1)}{2\log_{2}{\left(1-\frac{1}{2d}\right)}-(1-\epsilon_n)} \log_{2}{\left(1-\frac{1}{2d}\right)}}\\
    &\quad  \times \frac{1}{(\log{n})^{\frac{2(\epsilon_n-1)}{2\log_{2}{\left(1-\frac{1}{2d}\right)}-(1-\epsilon_n)} \frac{2\epsilon_n}{1-\epsilon_n}\log_{2}{\left(1-\frac{1}{2d}\right)}}}\\
&= c_{1} n^{\frac{2(\epsilon_n-1)}{2\bigl(\frac{-\frac{1}{2d}}{\log{2}}\bigr)-(1-\epsilon_n)} \bigl(\frac{-\frac{1}{2d}}{\log{2}}\bigr)}\\
&\quad \times \frac{1}{(\log{n})^{\frac{2(\epsilon_n-1)}{2\log_{2}\bigl(1-\frac{1}{2d}\bigr)-(1-\epsilon_n)} \frac{2\epsilon_n}{1-\epsilon_n}\log_{2}\bigl(1-\frac{1}{2d}\bigr)}}\\
&= c_{1}n^{-\bigl(\frac{1-\epsilon_n}{1+(1-\epsilon_n)d \log{2}}\bigr)} (\log{n})^{\bigl(\frac{\epsilon_n}{ 1+d \log{2}(1-\epsilon_n )}\bigr)}.
\end{align*}

and, for the second term, with the same arguments

\[ \frac{\tilde{c}}{\epsilon_n} n^{(c(d)-1)(1-\epsilon_n)} \log{n}^{\epsilon_n(1-c(d))} =  \frac{\tilde{c}}{\epsilon_n}    n^{-\bigl(\frac{1-\epsilon_n}{1+(1-\epsilon_n)d \log{2}}\bigr)}  (\log{n})^{\bigl(\frac{\epsilon_n}{ 1+d \log{2}(1-\epsilon_n )}\bigr)} \]\\
for a constant $\tilde{c}$ independent of $\epsilon_n,$
hence,

\[   \mathbb{E}(\tilde{m}^{Cen}_{\infty,n}(x)- m(x) )^2 \leq \frac{C}{\epsilon_n}  n^{-\bigl(\frac{1-\epsilon_n}{1+(1-\epsilon_n)d \log{2}}\bigr)}  (\log{n})^{\bigl(\frac{\epsilon_n}{ 1+d \log{2}(1-\epsilon_n )}\bigr)}  , \]
and consequently,
\begin{align*}
\frac{C}{\epsilon_n}  n^{-\bigl(\frac{1-\epsilon_n}{1+(1-\epsilon_n)d \log{2}}\bigr)}  (\log{n})^{\bigl(\frac{\epsilon_n}{ 1+d \log{2}(1-\epsilon_n )}\bigr)} &= \frac{C}{\epsilon_n} n^{-\bigl(\frac{1}{1+d \log{2}}\bigr)}(\log{n})^{\bigl(\frac{\epsilon_n}{ 1+d \log{2}(1-\epsilon_n )}\bigr)}\\
&\quad \times n^{\bigl(\frac{\epsilon_n}{ (1+d \log{2}) (1 + (1-\epsilon_n ))d \log{2} }\bigr)}\\
&\leq  \frac{C}{\epsilon_n} n^{-\bigl(\frac{1}{1+d \log{2}}\bigr)}(\log{n})^{\bigl(\frac{\epsilon_n}{d \log{2}(1-\epsilon_n )}\bigr)}\\
&\quad \times (\log{n})^{ \frac{\log{n}}{\log{\log{n}}}
 \bigl(\frac{\epsilon_n}{ (d \log{2})^{2} (1-\epsilon_n ) }\bigr)} .\\
\end{align*}
Finally we finish the proof by selecting $\epsilon_n= \frac{1}{\log{n}},$

and 
\[  \mathbb{E}(\tilde{m}^{Cen}_{\infty,n}(x)- m(x) )^2 \leq \tilde{C} n^{-\bigl(\frac{1}{1+d \log{2}}\bigr)} (\log{n} ). \]

\end{proof}
\end{theorem}

\begin{theorem}
  $   \textbf{Y}=m(\textbf{X}) +\epsilon$ where $\epsilon$ is a zero mean Gaussian noise with finite variance independent of $\textbf{X} $. Assuming also that $\textbf{X} $ is uniformly distributed in $ [0,1]^d$ and $m$ is a Lipschitz function. Providing $k\to \infty$, there exists a constant  $ \tilde{C}$  such that for every $n>1$ , for every $x \in [0,1]^d$ 

\[  \mathbb{E}(\tilde{m}^{Un}_{\infty,n}(x)- m(x) )^2 \leq \tilde{C} n^{-\bigl(\frac{1}{1+\frac{3}{2}d \log{2}}\bigr)} (\log{n} ). \]

    \begin{proof}
By arguing with the same reasoning as the proof of the centered random forest we can verify that 
\[  \bigg( \mathbb{P}(C_{t}^c(x))  \bigg)^{\frac{1}{2}}     \leq    \bigg( \mathbb{P}( \rvert A_{n}(x) \rvert >t)+  \mathbb{P}( \lvert B_{n}(x) \rvert >t)  \bigg)^{\frac{1}{2}} \leq c_{3} e^{- c_{4}\frac{nt^2}{ 2^k}}.\] 
 for some constants $c_{3}, c_{4}$ independent of $k$ and $n.$
 The rate of consistency for the Uniform KeRF is the minimum of the right hand in the inequality in terms of $n$

   \[ \mathbb{E}( \tilde{m}^{Un}_{\infty,n} -m(x) )^{2} \leq c_{1} \left(1-\frac{1}{3d}\right)^{2k} + c_{2}t^{2} + c_{3}\log{n} e^{-c_{4}t^{2}\frac{n}{2^k}}.\]

    We compute the minimum of the right-hand side of the inequality for $t \in (0,1),$ 
        \begin{align*}
            2c_{2}t - 2tc_{4}\log{n}c_{3} \frac{n}{2^k} e^{-c_{4}t^{2}\frac{n}{2^k}} &= 0 \quad \Rightarrow\\
            e^{-c_{4}t^{2}\frac{n}{2^k}} &= \frac{c_{2}}{c_{3}c_{4}} \frac{2^k}{n\log{n}} \quad \text{and}\\
            t^{2} &= \frac{1}{c_{4}} \frac{2^k}{n}\log{\left(\frac{c_{2}}{c_{3}c_{4}} \frac{n \log{n}}{2^k}\right)}.
        \end{align*}
        
        Hence, the inequality becomes,
        
        \begin{align*}
             \mathbb{E}(\tilde{m}^{Un}_{\infty,n}(x)- m(x) )^2 &\leq c_{1} \left(1-\frac{1}{3d}\right)^{2k} + c_{2} \frac{1}{c_{4}} \frac{2^k}{n}\log{\left(\frac{c_{2}}{c_{3}c_{4}} \frac{n \log{n}}{2^k}\right)}\\
             &\quad \quad + c_{3} \log{n}\frac{c_{2}}{c_{3}c_{4}}\frac{2^k}{n\log{n}}\\
             &= c_{1} \left(1-\frac{1}{3d}\right)^{2k} + c_{2} \frac{1}{c_{4}} \frac{2^k}{n}\log{\left(\frac{c_{2}}{c_{3}c_{4}} \frac{n \log{n}}{2^k} e^{\frac{c_{2}}{c_{4}}}\right)} .
        \end{align*}
        
 For every $\epsilon_n \in (0,2] $ it holds, $\log{x} \leq \frac{1}{\epsilon_n} x^{\epsilon_n}.$ Then one has that,

    \[\mathbb{E}( \tilde{m}^{Un}_{\infty,n} -m(x) )^{2} \leq c_{1} \left(1-\frac{1}{3d}\right)^{2k} +  \frac{c_2 (   e^{\frac{c_{2}}{c_{4}}}  \frac{c_{2}}{c_{3}c_{4}} )^{n} }{c_4 \epsilon_n}\left(\frac{2^k}{n} (\log{n})^{\frac{\epsilon_n}{1-\epsilon_n}}\right)^{1-\epsilon_n}.\]

    We pick $k =  c(d) \log_{2}{\frac{n}{(\log{n})^{\frac{\epsilon_n}{1-\epsilon_n}}}}$

    Therefore, 

\[     \frac{c_2 (   e^{\frac{c_{2}}{c_{4}}}  \frac{c_{2}}{c_{3}c_{4}} )^{n} }{c_4 \epsilon_n}\left(\frac{2^k}{n} (\log{n})^{\frac{\epsilon_n}{1-\epsilon_n}}\right)^{1-\epsilon_n} \leq \frac{c'}{\epsilon_n} n^{(c(d)-1)(1-\epsilon_n)} \log{n}^{\epsilon_n(1-c(d))} ,\]

for a constant $ c'$ independent of $n$ and,

\begin{align*}
 c_{1} \left(1-\frac{1}{3d}\right)^{2k} &= c_{1} \left(1-\frac{1}{3d}\right)^{2c(d) \log_{2}{\frac{n}{(\log{n})^{\frac{\epsilon_n}{1-\epsilon_n}}}} }\\
 &= c_{1} 2^{2c(d) \log_{2}{\left(1-\frac{1}{3d}\right)}\log_{2}{\frac{n}{(\log{n})^{\frac{\epsilon_n}{1-\epsilon_n}}}} }\\
 &= c_{1} n^{2c(d)  \log_{2}{\left(1-\frac{1}{3d}\right)}} \frac{1}{(\log{n})^{c(d) \frac{2\epsilon_n}{1-\epsilon_n}\log_{2}{\left(1-\frac{1}{3d}\right)}} }.\\
\end{align*}

Therefore,

\[ c(d)= \frac{\epsilon_n-1}{2\log_{2}{\left(1-\frac{1}{3d}\right)}-(1-\epsilon_n)}\]
 Finally,
\begin{align*}
     c_{1}n^{2c(d) \log_{2}{\left(1-\frac{1}{3d}\right)}} \frac{1}{(\log{n})^{c(d) \frac{2\epsilon_n}{1-\epsilon_n}\log_{2}{\left(1-\frac{1}{3d}\right)}}}
&= c_{1} n^{\frac{2(\epsilon_n-1)}{2\log_{2}{\left(1-\frac{1}{3d}\right)}-(1-\epsilon_n)} \log_{2}{\left(1-\frac{1}{3d}\right)}}\\
    &\quad  \times \frac{1}{(\log{n})^{\frac{2(\epsilon_n-1)}{2\log_{2}{\left(1-\frac{1}{3d}\right)}-(1-\epsilon_n)} \frac{2\epsilon_n}{1-\epsilon_n}\log_{2}{\left(1-\frac{1}{3d}\right)}}}\\
&= c_{1} n^{\frac{2(\epsilon_n-1)}{2\bigl(\frac{-\frac{1}{3d}}{\log{2}}\bigr)-(1-\epsilon_n)} \bigl(\frac{-\frac{1}{3d}}{\log{2}}\bigr)}\\
&\quad\times \frac{1}{(\log{n})^{\frac{2(\epsilon_n-1)}{2\left(\frac{-\frac{1}{3d}}{\log{2}}\right)-(1-\epsilon_n)} \frac{2\epsilon_n}{1-\epsilon_n}\left(\frac{-\frac{1}{3d}}{\log{2}}\right)}}\\
&= n^{-\left(\frac{2(1-\epsilon_n)}{1+(1-\epsilon_n)d \log{2}}\right)} \frac{1}{(\log{n})^{\frac{2\epsilon_n}{-2+3d\log{2}(\epsilon_n-1)}}} \\
&= n^{-\left(\frac{2(1-\epsilon_n)}{2+(1-\epsilon_n)3d \log{2}}\right)} (\log{n})^{\left(\frac{2\epsilon_n}{ 2+3d \log{2}(1-\epsilon_n)}\right).}
\end{align*}

and, for the second term, with the same arguments

\[ \frac{\tilde{c}}{\epsilon_n} n^{(c(d)-1)(1-\epsilon_n)} \log{n}^{\epsilon_n(1-c(d))} =  \frac{\tilde{c}}{\epsilon_n}    n^{-\bigl(\frac{1-\epsilon_n}{1+(1-\epsilon_n)\frac{3}{2}d \log{2}}\bigr)}  (\log{n})^{\bigl(\frac{\epsilon_n}{ 1+d \frac{3}{2}\log{2}(1-\epsilon_n )}\bigr)} , \]\\
for a constant $\tilde{c}$ independent of $\epsilon_n$
hence,

\[   \mathbb{E}(\tilde{m}^{Un}_{\infty,n}(x)- m(x) )^2 \leq \frac{C}{\epsilon_n}  n^{-\bigl(\frac{1-\epsilon_n}{1+(1-\epsilon_n)\frac{3}{2}d \log{2}}\bigr)}  (\log{n})^{\bigl(\frac{\epsilon_n}{ 1+\frac{3}{2}d \log{2}(1-\epsilon_n )}\bigr)}  , \]
and consequently,
\begin{align*}
\frac{C}{\epsilon_n}  n^{-\bigl(\frac{1-\epsilon_n}{1+(1-\epsilon_n)\frac{3}{2}d \log{2}}\bigr)}  (\log{n})^{\bigl(\frac{\epsilon_n}{ 1+\frac{3}{2}d \log{2}(1-\epsilon_n )}\bigr)} &= \frac{C}{\epsilon_n} n^{-\bigl(\frac{1}{1+\frac{3}{2}d \log{2}}\bigr)}(\log{n})^{\bigl(\frac{\epsilon_n}{ 1+\frac{3}{2}d \log{2}(1-\epsilon_n )}\bigr)}\\
&\quad \times n^{\bigl(\frac{\epsilon_n}{ (1+\frac{3}{2}d \log{2}) (1 + (1-\epsilon_n ))d \log{2} }\bigr)}\\
&\leq  \frac{C}{\epsilon_n} n^{-\bigl(\frac{1}{1+\frac{3}{2}d \log{2}}\bigr)}(\log{n})^{\bigl(\frac{\epsilon_n}{\frac{3}{2}d \log{2}(1-\epsilon_n )}\bigr)}\\
&\quad \times  (\log{n})^{ \frac{\log{n}}{\log{\log{n}}}
 \bigl(\frac{\epsilon_n}{ (\frac{3}{2}d \log{2})^{2} (1-\epsilon_n ) }\bigr)}.\\
\end{align*}
Finally we finish the proof by selecting $\epsilon_n= \frac{1}{\log{n}},$

and 
\[  \mathbb{E}(\tilde{m}^{Un}_{\infty,n}(x)- m(x) )^2 \leq \tilde{C} n^{-\bigl(\frac{1}{1+\frac{3}{2}d \log{2}}\bigr)} (\log{n} ) .\]

\end{proof}

\end{theorem}

\section{Plots and Experiments.}

In the following section, we summarize the rates of convergence for the centered KeRF and the uniform KeRF compared with the minimax rate of convergence over the class of the Lipschitz functions \cite{minimax}. In addition, we provide numerical simulations where we compare the $L_2-$ error for different choices of the tree depth. All experiments performed with the software Python \url{ https://www.python.org/,} where random sets uniformly distributed in  $[0,1]^d$ have been created, for various examples for the dimension $d$ and the function $Y.$ For every experiment the set was divided in a training set(80 \%) and a testing set (20 \%) and afterwards the $L_2-$error was computed. For the centered KeRF we compare three different values of tree depth as they were provided in \cite{S},\cite{Interp}, and in Theorem \ref{TheoPetronio}. Moreover, for the uniform KeRF, we again compare the same values of tree depth as they were derived from \cite{S} and Theorem \ref{TheoAmbrogio} nevertheless, it is not known if the uniform-KeRF algorithm converges when our estimator function interpolates the data set. Of course, in practice, since real data might violate the assumptions of the theorems, one should try cross-validation for tuning the parameter of the algorithms.\\
    Comparing the rates of consistency for centered KeRF and the depth of the corresponding trees:
\begin{itemize}
  \item   Scornet in \cite[Theorem 1]{S} rate of convergence: $n^{-( \frac{1}{dlog2+3} )} (\log{n})^2,$ and $k=\lceil \frac{1}{\log2 +\frac{3}{d}}\log{\frac{n}{\log{n}^{2}}} \rceil  $
  \item New rate of convergence :  \[ n^{-\bigl(\frac{1}{1+d \log{2}}\bigr)}(\log{n})  , \; \text{and} \; k=\lceil \frac{\frac{1}{\log{n} } -1}{2\log_{2}(1-\frac{1}{2d})-( 1-\frac{1}{\log{n}})}\log_{2}\frac{n}{(\log{n})^{\frac{\frac{1}{\log{n}}}{1-\frac{1}{\log{n}}}}} \rceil \]
  \item Minimax \cite{minimax}
rate of consistency over the class of Lipschitz functions: $ n^{\frac{-2}{d+2}} $
functions
\end{itemize}

Thus, theoretically, the improved rate of consistency is achieved when trees grow at a deeper level.  
  \begin{figure}[H]
    \centering
    \includegraphics[width=0.6\textwidth]{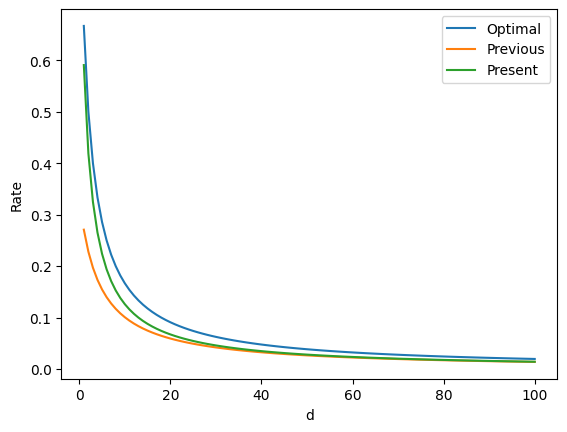}
    \caption{Plot of the exponents of n, for the previous rate of convergence for the centered KeRF algorithm, the new rate of convergence, and the optimal over the class of the Lipschitz functions.}
    \label{fig:plot_label}
\end{figure}

  Comparing the rates of convergence for uniform KeRF and the depth of the corresponding trees:
\begin{itemize}
\item Scornet in \cite[Theorem 2]{S}: rate of convergence : $n^{-( \frac{2}{3dlog2+6} )} (\log{n})^2 , $ and $k= \lceil \frac{1}{\log2 +\frac{3}{d}}\log{\frac{n}{\log{n}^{2}}} \rceil $
  \item New rate of convergence:   
 
$n^{-( \frac{2}{3d\log2+2} )}(\log{n}), $     and $  k=\lceil \frac{\frac{1}{\log{n} } -1}{2\log_{2}(1-\frac{1}{3d})-( 1-\frac{1}{\log{n}})}\log_{2}\frac{n}{(\log{n})^{\frac{\frac{1}{\log{n}}}{1-\frac{1}{\log{n}}}}} \rceil  $
  
  \item Minimax \cite{minimax}
rate of convergence for the consistency over the class of Lipschitz functions: $ n^{\frac{-2}{d+2}} $
functions
\end{itemize}
  
Thus, theoretically, as in the case of centered random KeRF the improved rate of consistency is achieved when trees grow at a deeper level.  

  \begin{figure}[H]
    \centering
    \includegraphics[width=0.6\textwidth]{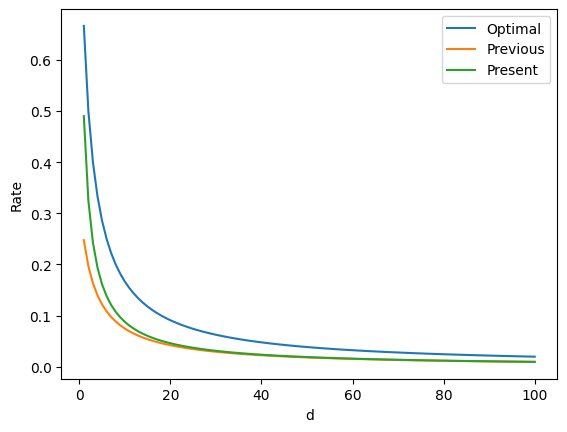}
     \caption{Plot of the exponents of n, for the previous rate of convergence for the uniform KeRF algorithm, the new rate of convergence, and the optimal over the class of the Lipschitz functions.}
    \label{fig:plot_label}
\end{figure}

 \begin{figure}[H]
    \centering
    \includegraphics[width=0.6\textwidth]{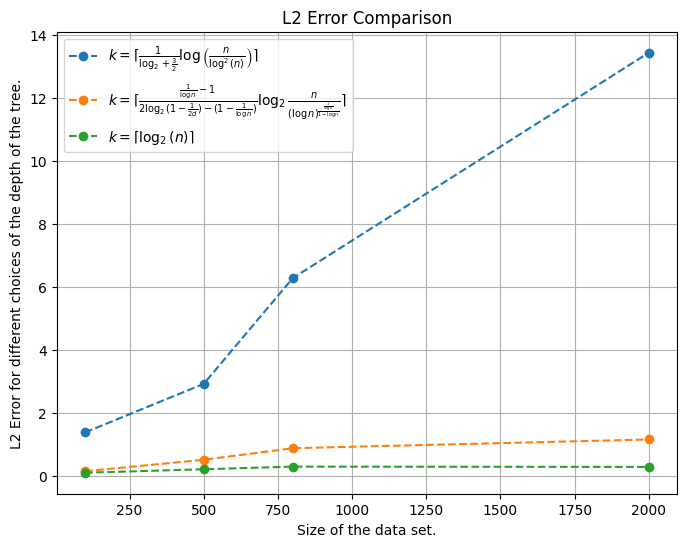}
     \caption{Plot of the $L_{2}-$error of the centered KeRF-approximation for three different values of $k$ for the function $Y=X^{2}_{1}+e^{-X^{2}_{2}} + \epsilon, $ where $\epsilon \sim \mathcal{N}(0,\frac{1}{2}),$ against different data set size.}
    \label{fig:plot_label}
\end{figure}

 \begin{figure}[H]
    \centering
    \includegraphics[width=0.6\textwidth]{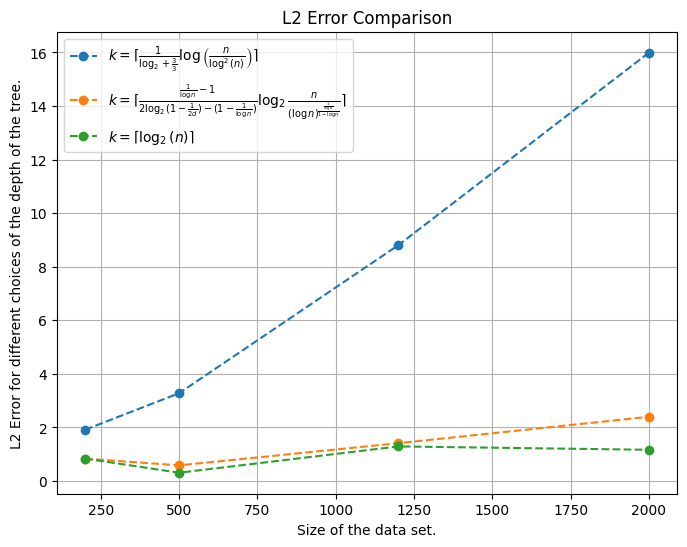}
     \caption{Plot of the $L_{2}-$error of the centered KeRF-approximation for three different values of $k$ for the function $Y=X^{2}_{1}+\frac{1}{e^{X^{2}_{2}} +e^{X^{2}_{3}}}+\epsilon $ where $ \epsilon \sim \mathcal{N}(0,0.5,)$ against different data set size.}
    \label{fig:plot_label}
\end{figure}

 \begin{figure}[H]
    \centering
    \includegraphics[width=0.6\textwidth]{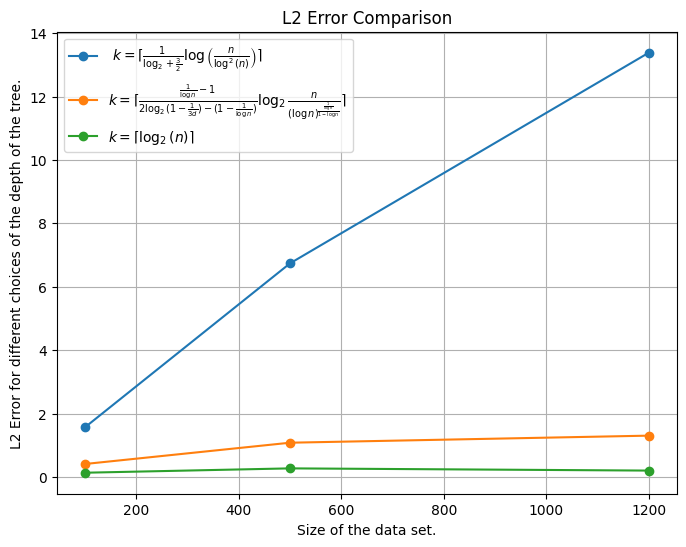}
     \caption{Plot of the $L_{2}-$error of the uniform KeRF-approximation for three different values of $k$ for the function $Y=X^{2}_{1}+e^{-X^{2}_{2}}+\epsilon$  where $ \epsilon \sim \mathcal{N}(0,0.5),$ against different data set size.}
    \label{fig:plot_label}
\end{figure}

 \begin{figure}[H]
    \centering
    \includegraphics[width=0.6\textwidth]{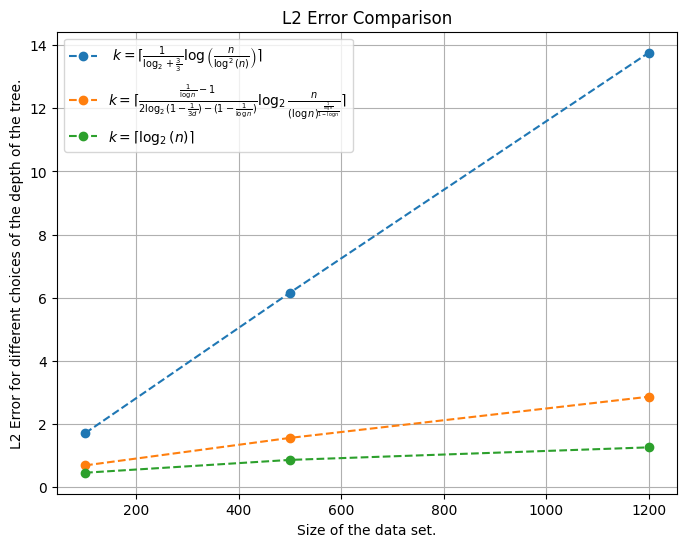}
     \caption{Plot of the $L_{2}-$error of the uniform KeRF-approximation for three different values of $k$ for the function $Y=X^{2}_{1}+\frac{1}{(e^{X^{2}_{3}} + e^{X^{2}_{2}}    )} + \epsilon $ where $ \epsilon \sim \mathcal{N}(0,0.5),$ against different data set size.}
    \label{fig:plot_label}
\end{figure}

\section{Analysis of the Kernel.}

\subsection{Fourier transforms on finite groups and related RKHS}
  Following the notation of \cite{Rudin} we recall here some basic notions from the general theory for a finite, abelian group $G$, endowed with its counting measure. The {\it dual group} $\Gamma=\widehat{G}$ of $G$ is populated by labels $a$ for homomorphisms $\gamma_a:G\to\mathbb{T}=\{e^{i t}:\ t\in\mathbb{R}\}$. Given a function $f:G\to\mathbb{C}$, its
{\it Fourier transform} $\widehat{f}:\Gamma\to\mathbb{C}$ is defined as
\begin{equation}
\widehat{f}(a)=\sum_{x\in G}f(x)\overline{\gamma_a(x)}.    
\end{equation}
We make $\Gamma$ into a (finite), additive group by setting
\[
\gamma_{a+b}=\gamma_a\cdot\gamma_b,\text{ and }\gamma_x(a):=\gamma_a(x).
\]
It turns out that they have the same number of elements, $\sharp(G)=\sharp(\Gamma)$. Some basic properties are:
\begin{eqnarray*}
f(x)&=&\frac{1}{\sharp(\Gamma)}\sum_{a\in \Gamma}\widehat{f}(a)\gamma_a(x)\text{ (inverse Fourier transform),}\crcr  
\sum_{x\in G}|f(x)|^2&=&\frac{1}{\sharp(\Gamma)}\sum_{a\in \Gamma}|\widehat{f}(a)|^2\text{ (Plancherel),}\crcr 
\widehat{f\ast g}&=&\widehat{f}\cdot\widehat{g},
\end{eqnarray*}
where
\begin{equation}\label{eqConvOne}
    (f\ast g)(x)=\sum_{y\in G}f(x-y)g(y).
\end{equation}
We write
\begin{equation}\label{eqCheck}
    \check{\varphi}(x)=\sharp(\Gamma)^{-1}\sum_{a\in\Gamma}\varphi(a)\gamma_a(x),\text{ so that }\widehat{\check{\varphi}}=\varphi.
\end{equation}
The unit element of convolution in $G$ is $\delta_0$.

In the other direction, for $\varphi,\psi:\Gamma\to\mathbb{C}$ we define
\begin{equation}\label{eqConvTwo}
    (\varphi\ast\psi)(a)=\frac{1}{\sharp(\Gamma)}\sum_{b\in\Gamma}\varphi(a-b)\psi(b),
\end{equation}
and similarly to above, $\widehat{\check{\varphi}\check{\psi}}=\varphi\ast\psi$. The unit element on convolution in $\Gamma$
is $\sharp(\Gamma)\delta_0$.

A function $\varphi$ on $\Gamma$ is {\it positive definite} if
\[
\sum_{a,b\in\Gamma}^n c(a) \overline{c(b)}\varphi(b-a)\ge0.
\]
\begin{theorem}\label{TheoBochner}[Bochner's Theorem]
    A function $\varphi:\Gamma\to\mathbb{C}$ is positive definite if and only if there exists $\mu:G\to\mathbb{R}_+$ such that 
    $\varphi=\widehat{\mu}$.
\end{theorem}
The proof for the finite group case is elementary, but we include it because it highlights the relationship between the measure $\mu$ and the functions $\varphi$.
\begin{proof}
[{\bf If.}]
\begin{eqnarray}
   \sharp(\Gamma)^{-2}\sum_{a,b\in\Gamma}\widehat{\mu}(b-a) c(a)\overline{c(b)}&=&\sum_{x\in G}\sharp(\Gamma)^{-2}\sum_{a,b\in\Gamma}\mu(x)\overline{\gamma_{b-a}(x)}c(a)\overline{c(b)}\crcr
   &=&\sum_{x\in G}\sharp(\Gamma)^{-2}\sum_{a,b\in\Gamma}\mu(x)\overline{\gamma_{b}(x)}\overline{c(b)}\gamma_a(x)c(a)\crcr
   &=&\sum_{x\in G}\mu(x)\left|\sharp(\Gamma)^{-1}\sum_{a\in\Gamma}c(a)\gamma_a(x)\right|^2\crcr
   &=&\sum_{x\in G}\mu(x)\left|\check{c}(x)\right|^2\ge0.
\end{eqnarray}
{\bf Only if.} Since for all $b$ in $\Gamma$,
\begin{eqnarray}\label{eqInvBoch}
\mu(x)\sharp(\Gamma)&=&\sum_{a\in\Gamma}\varphi(a)\gamma_x(a)=\sum_{a\in\Gamma}\varphi(a-b)\gamma_x(a-b)\crcr
&=&\sum_{a\in\Gamma}\varphi(a-b)\gamma_x(a)\overline{\gamma_x(b)},
\end{eqnarray}
we have
\begin{equation}\label{eqInvBochTwo}
    \mu(x)\sharp(\Gamma)^2=\sum_{a,b\in\Gamma}\varphi(a-b)\gamma_x(a)\overline{\gamma_x(b)}\ge0,
\end{equation}
by the assumption.
\end{proof}
We now come to reproducing kernels on $\Gamma$ which are based on positive definite functions $\varphi:\Gamma\to\mathbb{R}_+$.
Set
\begin{equation}\label{eqK}
    K(a,b)=\varphi(a-b)=K_b(a),\ K:\Gamma\times\Gamma\to\mathbb{C},
\end{equation}
and set
\begin{equation}\label{eqHK}
H_K=\text{span}\{K_b:\ b\in\Gamma\}\ni\sum_{b\in\Gamma}c(b)K_b,
\end{equation}
where $H_K$ is the Hilbert space having $K$ as reproducing kernel. We wish to have a more precise understanding of it.

We start by expressing the norm of an element on $H_K$ is several equivalent ways,
\begin{eqnarray}\label{eqNormHK}
    \left\|\sum_{b\in\Gamma}c(b)K_b\right\|_{H_K}^2&=&\sum_{a,b\in\Gamma}\overline{c(a)}c(b)\langle K_b,K_a\rangle\crcr 
    &=&\sum_{a,b\in\Gamma}\overline{c(a)}c(b)K(a,b)=\sum_{a,b\in\Gamma}\overline{c(a)}c(b)\widehat{\mu}(a-b)\crcr
    &=&\sum_{a,b\in\Gamma}\overline{c(a)}c(b)\sum_{x\in G}\mu(x)\gamma_{b-a}(x)\crcr
    &=&\sum_{x\in G}\mu(x)\sum_{a,b\in\Gamma}\overline{c(a)}c(b)\gamma_{b}(x)\overline{\gamma_a(x)}\crcr
    &=&\sum_{x\in G}\mu(x)\left|\sum_{b\in\Gamma}c(b)\gamma_b(x)\right|^2\crcr
    &=&\sharp(\Gamma)^2\sum_{x\in G}\mu(x)\left|\check{c}(x)\right|^2=\sharp(\Gamma)^2\sum_{x\in G}\left|\mu(x)^{1/2}\check{c}(x)\right|^2.
\end{eqnarray}
In other terms,
\begin{equation}\label{eqIsometry}
    \sharp(\Gamma)^{-1}\sum_{b\in\Gamma}c(b)K_b\mapsto \check{c}
\end{equation}
is an isometry of $H_K$ onto $L^2(\mu)$. This will become important later, when we verify that for our kernels $\text{supp}(\mu)$ is sparse in $G$. In fact, $\text{dim}(H_K)=\sharp(\text{supp}(\mu))$.
\begin{corollary}\label{corAllThat}
    As a linear space, $H_K$ is determined by $\text{supp}(\mu)$:
    \[
    \psi\in H_K\text{ if and only if }\text{supp}(\check{\psi})\subseteq\text{supp}(\mu).
    \]
\end{corollary}
Let $E\subseteq G$. We denote
\begin{equation}\label{eqLinear}
    L_E=\{G\xrightarrow{\psi}\mathbb{C}:\ \text{supp}(\check{\psi})\subseteq E\}.
\end{equation}

Next, we look for the natural orthonormal system provided by the Fourier isometry \eqref{eqIsometry}. Fr $x\in G$, let $\check{c_x}=\mu(x)^{-1/2}\delta_x$: $\{\check{c_x}:\ x\in E:=\text{supp}(\mu)\}$ is a orthonormal system for $L^2(\mu)$,
and so $\{e_x:\ x\in E\}$ is an orthonormal basis for $H_K$, where
\begin{eqnarray}\label{eqPrima}
    c_x(b)=\sum_{y\in G}\mu(x)^{-1/2}\delta_x(y)\overline{\gamma_b(y)}=\mu(x)^{-1/2}\overline{\gamma_b(x)},
\end{eqnarray}
and
\begin{eqnarray}\label{eqSeconda}
    e_x(a)&=&\sharp(\Gamma)^{-1}\sum_{b\in\Gamma}c_x(b)K_b(a)\crcr 
    &=&\frac{\mu(x)^{-1/2}}{\sharp(\Gamma)}\sum_{b\in\Gamma}K_b(a)\overline{\gamma_b(x)}\crcr
    &=&\frac{\mu(x)^{-1/2}}{\sharp(\Gamma)}\sum_{b\in\Gamma}\varphi(a-b)\overline{\gamma_b(x)}\crcr
    &=&\frac{\mu(x)^{-1/2}}{\sharp(\Gamma)}\sum_{b\in\Gamma}\varphi(a-b)\gamma_{a-b}(x)\overline{\gamma_a(x)}\crcr
    &=&\mu(x)^{-1/2}\mu(x)\overline{\gamma_a(x)}\crcr
    &=&\mu(x)^{1/2}\overline{\gamma_a(x)}.
\end{eqnarray}
Let's verify that we obtain the reproducing kernel from the o.n.b. as expected,
\begin{eqnarray}\label{eqOnbK}
\sum_{x\in\Gamma}e_x(a)\overline{e_x(b)}&=&\sum_{x\in\Gamma}\mu(x)\overline{\gamma_x(a)}\gamma_x(b)\crcr 
&=&\sum_{x\in\Gamma}\mu(x)\overline{\gamma_x(a-b)}\crcr 
&=&\widehat{\mu}(a-b)\crcr
&=&\varphi(a-b).
\end{eqnarray}
\begin{remark}
    Since any finite, abelian group can be written as the direct product of cyclic groups,
    \begin{equation}\label{eqCyclic}
        G=\bigoplus_{l=1}^L \mathbb{Z}_{m_l},
    \end{equation}
    its dual $\Gamma$ can be written in the same way, because $\widehat{\mathbb{Z}_m}\equiv\mathbb{Z}_m$. From the Fourier point of view, the only difference is that, if on $G$ we consider the counting measure, then on $\Gamma$ we consider normalized counting measure, as we did above.
\end{remark}
\subsection{A reproducing kernel from the Centred KeRF}

\subsubsection{The Fourier analysis of the kernel}
We identify every real number $x \in [0,1]$ with its dyadic decomposition $x=0.x_{1}x_{2}x_{3}...$ with $x_{i} \in \{ 0,1\}$ for $i \in \mathbb{N}.$

Here we consider the group 
\begin{equation}\label{eqZnp}
G=\mathbb{Z}_2^{kd}\ni x=(x_i^j)_{\ontopof{i=1,\dots,k}{j=1,\dots,d}}=(x^1|x^2|\dots|x^d)=\begin{pmatrix}
x_1\crcr
x_2\crcr
\dots\crcr
x_k
\end{pmatrix}.  
\end{equation}
The kernel $K:\Gamma\times\Gamma\to\mathbb{C}$ corresponding to the kernel $K^{cen}_{k} $ is,
\begin{eqnarray}\label{eqOurKernel}
  K(a,b)&=&\sum_{\ontopof{l\in\mathbb{N}^d}{|l|=k}}\frac{1}{d^k}\binom{k}{l}\prod_{j=1}^d\chi\left(a^j_1=b^j_1,\dots,a^j_{k_j}=b^j_{k_j}\right)\crcr
  &=&\sum_{\ontopof{l\in\mathbb{N}^d}{|l|=k}}\frac{1}{d^k}\binom{k}{l}\prod_{j=1}^d\prod_{i=1}^{k_j}\chi\left(a^j_i=b^j_i\right)\crcr
&=&\varphi(a-b), 
\end{eqnarray}
where $\binom{k}{l} $ is the multidimensional binomial coefficient and $\chi $  the characteristic function. For the last equality, we consider the dyadic representation of a number in $(0,1]$ whose digits are not definitely vanishing. The fact that $0$ does not have such representation is irrelevant since the probability that one of the coordinates of the data vanishes is zero.

We now compute the anti-Fourier transform $\mu=\check{\varphi}$. We know that $\sharp(\Gamma)=2^{kd}$, and that the characters of $\mathbb{Z}_2^{kd}$ have the form
\begin{equation}\label{eqCharZtwo}
    \gamma_a(x),\ x\in \mathbb{Z}_2^{kd},\ a\in \widehat{\mathbb{Z}_2^{kd}},\ a\cdot x=a^1_1 x^1_1+\dots+a^d_k x^d_k.
\end{equation}
Hence,
\begin{eqnarray}\label{eqAFTvarphi}
2^{kd}p^n\mu(x)&=&d^k\sum_{a\in\Gamma}\varphi(a)\gamma_a(x)\crcr 
&=&d^k\sum_{a\in\Gamma}\varphi(a)(-1)^{a\cdot x}\crcr 
&=&\sum_{a\in\Gamma}\sum_{\ontopof{l\in\mathbb{N}^d}{|l|=k}}\binom{k}{l}\prod_{j=1}^d\prod_{i=1}^{k_j}\chi\left(a^j_i=0\right) (-1)^{a\cdot x}\crcr 
&=&\sum_{a\in\Gamma}\sum_{\ontopof{l\in\mathbb{N}^d}{|l|=k}}\binom{k}{l}\prod_{j=1}^d(-1)^{\tilde{a}_j^{k_j}\cdot\tilde{x}_j^{k_j}}\prod_{i=1}^{k_j}\left[\chi\left(a^j_i=0\right)(-1)^{a^j_i x^j_i}\right]\crcr
&\ &\text{where }\tilde{a}_j^{k_j}=\begin{pmatrix}
  a^j_{k_j+1}\crcr
  \dots\crcr 
  a^j_n
\end{pmatrix}\text{ is the lower, $(k-k_j)$-dimensional} \crcr
&\ &\text{part of the column $a^j$,}\crcr
&=&\sum_{\ontopof{l\in\mathbb{N}^d}{|l|=k}}\binom{k}{l}\prod_{j=1}^d(-1)^{\tilde{a}_j^{k_j}\cdot\tilde{x}_j^{k_j}}\prod_{i=1}^{k_j}\chi\left(a^j_i=0\right)\crcr
&=&\sum_{\ontopof{l\in\mathbb{N}^d}{|l|=k}}\binom{k}{l}\sum_{\ontopof{a\in\Gamma}{a^1_1=\dots a^1_{k_1}=a^2_1=\dots=a^d_{k_d}=0}}\prod_{j=1}^d(-1)^{\tilde{a}_j^{k_j}\cdot\tilde{x}_j^{k_j}}.
\end{eqnarray}
The last expression vanishes exactly when for all $l$, there are some $1\le j\le d$, and some $k_j+1\le i\le k$ such that $x^j_i=1$,
due to the presence of the factor $(-1)^{a^j_i x^j_i}=(-1)^{a^j_i}$ which takes values $\pm1$ on summands having, two by two, the same absolute values.

If, on the contrary, there is $l$ such that for all $1\le j\le d$, and $k_j+1\le i\le k$, we have that $x^j_i=0$, then $\mu(x)\ne0$.
Since $|l|=k$ and there are $kd$ binary digits involved in the expression of $x$, the latter occurs exactly when the binary matrix representing $x$ has a large lower region in which all entries are $0$. More precisely, the number of vanishing entries must be at least
\begin{equation}\label{eqBinLower}
(k-k_1)+\dots+(k-k_p)=(d-1)k.    
\end{equation}
The number $N(d,k)$ of such matrices is the dimension of $H_K$, the Hilbert space having $K$ as a reproducing kernel. 

Next, we prove some estimates for the dimension of the reproducing kernel Hilbert space.

We summarize the main items in the following statement.
\begin{theorem}\label{TheoKFourier}
    Let $K:\Gamma\times\Gamma\to\mathbb{C}$ be the kernel in \eqref{eqOurKernel}, $K(a,b)=\varphi(a-b)$, and let
    \begin{equation}\label{eqKSupport}
    E_K=\text{supp}(\check{\varphi})\in K.
    \end{equation}
    Then,
\begin{enumerate}
    \item[(i)] as a linear space, $H_K=L_{E_K}$, where 
    \begin{eqnarray}\label{eqSuppK}
        E_K&=&\{x=(x^1|\dots|x^d):\ x^j_i=0 \text{ for }k_j+1\le i\le k,\text{ for some }l\crcr 
        &=&(k_1,\dots,k_d)\in\mathbb{N}^d \text{ with }k_1+\dots+k_d=k\};
    \end{eqnarray}
    \item[(ii)] For $x\in E_K$,
    \begin{equation}\label{eqMuHat}
        \check{\varphi}(x)=\frac{1}{2^{k}d^k}\sum_{\ontopof{l\in\mathbb{N}^d,\ |l|=k}{x^j_i=0 \text{ for }k_j+1\le i\le k}}\binom{k}{l}
    \end{equation}
\end{enumerate}    
\end{theorem}
To obtain the expression on \eqref{eqMuHat}, we used the fact that 
\[
\sharp\{a:\ a^1_1=\dots a^1_{k_1}=a^2_1=\dots=a^p_{k_p}=0\}=2^{(d-1)k}.
\]

\subsection{Some properties of $H_K$.}
\paragraph{\bf Linear relations.} 
Among all functions $\psi:\Gamma\to\mathbb{C}$, those belonging to $H_K$ (i.e., those belonging to $L_{E_K}$) are characterized by a set of linear equations,
\begin{equation}\label{eqLEK}
0=2^{np}p^n\mu(x)=\sum_{\ontopof{k\in\mathbb{N}^p,\ |k|=n}{x^j_i=0 \text{ for }k_j+1\le i\le n}}\binom{n}{k}\text{ for }x\notin E_K.
\end{equation}

\paragraph{\bf Multipliers.} A {\it multiplier} of $H_K$ is a function $m:\Gamma\to\mathbb{C}$ such that $m\psi\in H_K$ whenever $\psi\in H_K$. 
\begin{proposition}\label{PropMultip}
  The space $H_K$ has no nonconstant multiplier.   
\end{proposition}
In particular, it does not enjoy the {\it complete Pick property}, which has been subject of intensive research for the past twenty-five years \cite{AMc}. 
\begin{proof}
The space $H_K$ coincides as a linear space with $L_{E_K}$. Let $\Lambda_E=\check{L_E}$, which is spanned by $\{\delta_x:\ x\in E\}$. 
Observe that, since $0=(0|\dots|0)\in E_K$, the constant functions belong to $H_K$, hence, any multiplier $m$ of $H_K$ lies in $H_K$; $m=m\cdot 1\in H_K$. 

Suppose $m$ is not constant. Then, $\check{m}(a)\ne0$ for some $a\in E_K$, $a\ne0$. 
    Let $a$ be an element in $E_K$ such that $\check{m}(a)\ne0$.
    Since $H_K\ni m\cdot\widehat{\delta_x}$ for all $x$ in $E_K$, and $m\cdot\widehat{\delta_x}=\widehat{\check{m}\ast\delta_x}$, we have that the support of $\check{m}\ast\delta_x$ lies in $H_K$. Now, $\check{m}\ast\delta_x(y)=\check{m}(x-y)$, hence, we have that,
    for any $x$ in $E_K$, $y=x-a$ lies in $E_K$ as well. This forces $a=0$, hence $m$ to be constant.  
\end{proof}

\subsection{Bounds for dimension and generating functions.}
\begin{theorem}
We have the estimates:
    \begin{equation}\label{eqPfixed}
\text{dim}(H_K)\sim\frac{2^{k-d+1}k^{d-1}}{(d-1)!},\text{ hence }\frac{\text{dim}(H_K)}{2^{kd}}\sim\frac{k^{d-1}}{2^{k-1}(d-1)!2^{k(d-1)}}.    
\end{equation}
\end{theorem}
\begin{proof}
Let $l_{1}, l_{2},...,l_{d}$ such that
\[0\leq l_{1}+l_{2}+...+l_{d}=m \leq k  \]
where $m $ is a parameter and let also $\lambda = \lvert j : l_{j} \geq 1 \rvert = \lvert \{ \text{stop 1-digits} \} \rvert = \lvert \{ \text{back-entries} \} \vert $ where $ \lvert \cdot \rvert$ denotes the size of the sets,
and of course we have that $0 \leq m \leq k $ and $0\leq \lambda \leq d,m.$ Goal to obtain a bound for 
\[ N(k,d)= \sum_{m=0}^k \sum_{\lambda=0}^{d \wedge m} 2^{m-\lambda} { d \choose \lambda} \lvert  \{ (l_{1},l_{2},...,l_{d}): l_{1} +l_{2}+...+l_{d}=m \rvert \quad\text{and} \quad \lvert \{ j : l_{j}=1 \}=1 \rvert \} .\]
Let $A(m,\lambda) $ the m-th coefficient of $x$  in the polynomial
\begin{align*}
(x+x^{2}+...x^{m}+...)^{\lambda} &=(x(1+x+x^{2}+...)^{\lambda}\\
 &= (x^{\lambda}(1+x+...)^{\lambda})\\
 &=\frac{x^{\lambda}}{(1-x)^{\lambda}}
\end{align*}
and $2^{m} A(m,\lambda) $ is the m-th coefficient of $x,$  for the fraction $\frac{(2x)^{\lambda}}{(1-2x)^{\lambda}},$
therefore $2^{m-\lambda} A(m,\lambda) $ is the m-th coefficient of $\frac{x^{\lambda}}{(1-2x)^{\lambda}}. $
Let's see the first sum, \\
$B(m,d)$ is the m-th coefficient of $x$ :
\begin{align*}
  \sum_{\lambda=0}^{d \wedge m }{ d \choose \lambda} 2^{m-\lambda } A(m,\lambda) &= \sum_{\lambda=0}^{d \wedge m}{ d \choose \lambda} ( \frac{x}{1-2x} )^{\lambda}\\
 &= ( 1+ \frac{x}{1-2x})^d\\
 &=(\frac{1-x}{1-2x})^d
\end{align*}
Again by the same combinatoric argument we are looking the $k$-th coefficient of the function
\[ f(x)=\frac{1}{1-x}(\frac{1-x}{1-2x})^d . \]

Back to the estimate,\\
Let $a_{k}$ the coefficient of the power series centered at $z=0.$
\[ \max_{ |z|=r } | f(z)| = \max_{ |z|=r } \bigg| \frac{(1-z)^{d-1}}{(1-2z)^{d}} \bigg| =\max_{ |z|=r }\bigg| \frac{1}{1-z} \bigg(\frac{1-z}{1-2z}\bigg)^{d} \bigg| \leq 2 \max_{ \theta \in ( -\pi, \pi) } \bigg| \frac{1-re^{i \theta}}{1-2re^{i \theta}}    \bigg| ^{d} \]
After some calculations since $r$ is fixed one has that the maximum is achieved for $\theta=0.$ 
So $ \max_{ |z|=r } | f(z)| \leq 2 (\frac{1-r}{1-2r})^{d}  $
Our estimation finally becomes : 

\begin{align*}
    |a_{k}| &\leq \frac{2  \big( \frac{1-r}{1-2r} \big)^{d}  }{r^{k}}\\
    &=\frac{2(1-r)^{d}}{r^{k} (1-2r)^{k} }\\
    &=k^{d} 2^{k}(\frac{1}{2} + \frac{1}{2k})^{d}, \quad \quad  ( \text{since,} \quad r=\frac{1}{2}( 1-\frac{1}{k})  )\\
    &  =k^{d} (1+ \frac{1}{k})^{d} 2^{k-d}.
\end{align*}
Thus an estimate for the dimension of $ H_K$ is
\[ \frac{|a_k|}{2^{kd}}\lesssim \frac{k^{d} (1+ \frac{1}{k})^{d} 2^{k(1-d)}}{2^d}\]
Another estimate about the dimension of $H_K$.
For $f(z)=\sum_{n=0}^\infty a_n z^n$ we have
\[
|a_n|\le\frac{\max\{|f\left(r e^{i t}\right)|:\ |t|\le\pi\}}{r^n}.
\]
Consider the function 
\[
f(z)=\frac{(1-z)^{d-1}}{(1-2z)^d}
\]
in $|z|<1/2$ and let $r=\frac{1-1/k}{2}$. Then,
\begin{eqnarray*}
|a_k|&\le&\frac{(3/2)^{d-1}}{(1/k)^d(1-1/k)^k2^{-k}}\crcr 
&\le&(3/2)^{d-1}2^ke k^d.
\end{eqnarray*}
Thus,
\[
\frac{|a_k|}{2^{kd}}\lesssim\frac{k^d (3/2)^d}{2^{k(d-1)}}.
\]
Recursively working out the generating function one gets the estimates in \eqref{eqPfixed}.
\end{proof}

\bibliography{mybibl}

\bibliographystyle{abbrv}

\end{document}